\documentclass[9pt,twocolumn]{IEEEtran}
\IEEEoverridecommandlockouts
 \usepackage[T1]{fontenc}
\usepackage{amssymb,color}

 \renewcommand{\baselinestretch}{1} 
 
\usepackage{graphicx,amsmath,amssymb,color}

\usepackage{epsfig}
\newcommand{\EQ}{\begin{eqnarray}}
\newcommand{\EN}{\end{eqnarray}}
\newcommand{\EQQ}{\begin{eqnarray*}}
\newcommand{\ENN}{\end{eqnarray*}}
\newcommand{\R}{\mathbb R}
\newcommand{\N}{{\cal N}}
\newcommand{\n}{{\rm n}}

\renewcommand{\O}{{\cal O}}

\newcommand{\col}{\operatorname{col}}

\newcommand{\bproof}{\noindent {\it Proof:} }
 
\newcommand{\eproof}{\hfill $\square$}

\newcommand{\br}{{\bf r}}

 \newcommand{\bone}{\boldsymbol{1}}

\usepackage[outdir=./]{epstopdf}
\graphicspath{{../eps/}}
\DeclareGraphicsExtensions{.eps}

\renewcommand{\t}{^{\mbox{\tiny\sf T}}}

\newcommand{\bremark}
{\medskip\begin{remark}
\begin{rm}}
\newcommand{\eremark}{ \end{rm}\hfill \rule{0mm}{0mm}
\end{remark} }
\newcommand{\btheorem}{\medskip\begin{theorem} \begin{it}}
\newcommand{\etheorem}{\end{it} 
\end{theorem} }
\newcommand{\blemma}{\medskip\begin{lemma} \begin{it} }
\newcommand{\elemma}{ \end{it} 
\end{lemma} }
\newcommand{\bcorollary}{\medskip\begin{corollary} \begin{it} }
\newcommand{\ecorollary}{ \end{it} \hfill\rule{1mm}{2mm}
\end{corollary} }
\newcommand{\bdefinition}{\medskip\begin{definition} }
\newcommand{\edefinition}{ 
\end{definition} }
\newcommand{\bproposition}{\medskip\begin{proposition} }
\newcommand{\eproposition}{\hfill \rule{1mm}{2mm}
\end{proposition} }
\newcommand{\bexample}{\medskip\begin{example} \begin{rm}}
\newcommand{\eexample}{ \end{rm} \hfill\rule{1mm}{2mm}
\end{example} }

\newcommand{\basm}{\medskip\begin{assumption} \begin{rm} }
\newcommand{\easm}{ \end{rm} \hfill\rule{1mm}{2mm} \medskip
\end{assumption} }

\newtheorem{theorem}{Theorem}[section]
\newtheorem{lemma}{Lemma}[section]

\newtheorem{definition}{Definition}[section]
\newtheorem{remark}{ Remark}[section]
\newtheorem{corollary}{ Corollary}[section]
\newtheorem{proposition}{ Proposition}[section]
\newtheorem{example}{Example}[section]
\newtheorem{assumption}{ Assumption}

\begin{document}

\title{Synchronization of Frequency Modulated Multi-Agent Systems}

\author{{Zhiyong Chen}
\thanks{The author is with the School of  Engineering,
        The University of Newcastle, Callaghan, NSW 2308, Australia.
        Tel: +61 2 4921 6352, Fax: +61 2 4921 6993
      Email: {\tt\small  zhiyong.chen@newcastle.edu.au}.
        }}

\date{}

\maketitle
 
 \begin{abstract}
Oscillation synchronization phenomenon is widely observed in natural systems through frequency modulated signals, especially in biological neural networks. Frequency modulation is also one of most widely used technologies in engineering. However, due to the technical difficulty, oscillations are always simplified as unmodulated sinusoidal-like waves in studying the synchronization mechanism in the bulky literature over the decades. It lacks mathematical principles, especially systems and control theories, for frequency modulated multi-agent systems. This paper aims to bring a new formulation of synchronization of frequency modulated multi-agent systems.  It develops new tools to solve the synchronization problem by addressing three issues including frequency observation in nonlinear frequency modulated oscillators subject to network influence, frequency consensus via network interaction subject to observation error, and a well placed small gain condition among them. The architecture of the paper consists of a novel problem formulation, rigorous theoretical development, and numerical verification.

 \end{abstract}


\section{Introduction}

\begin{figure}[h]
 \center  \includegraphics [width=\hsize]{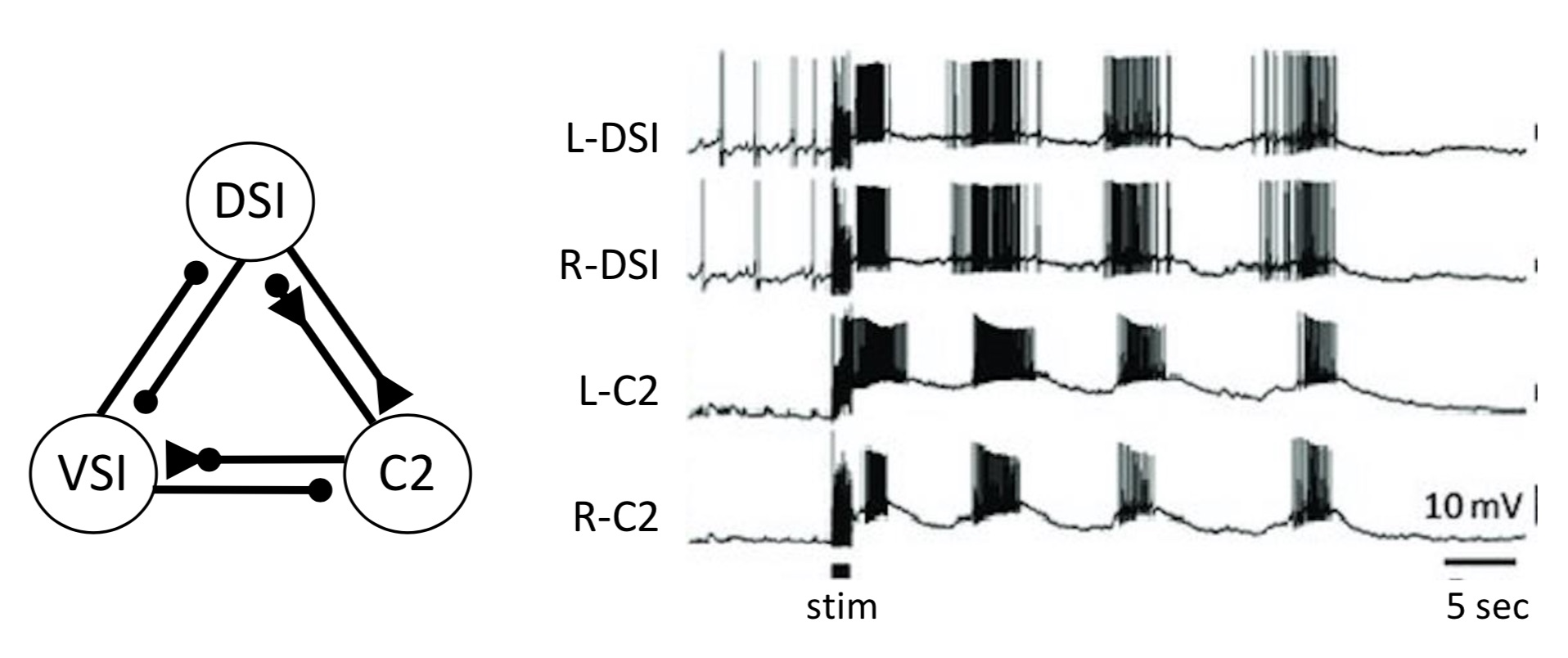}
\caption{ A real frequency modulated CPG model. The neural circuit and swim motor patters for the dorsal-ventral swimmer Tritonia. Left: The Tritonia swim CPG consists of three neuron types: DSI, C2, and VSI. Right: Simultaneous intracellular micro-electrode recordings show that two contralateral DSIs fire bursts of APs in phase with each other and slightly ahead of the two C2s. VSI (not recorded here) fires APs in the inter-burst interval.   
Courtesy of PNAS \cite{newcomb2012}.
     } \label{fig:slug}
\end{figure}

Frequency modulated signals extensively exist in natural network systems such as neural circuits. 
Neuroscientists are interested in modeling nervous systems from simple neural networks to 
even the human brain.  
One of the basic areas is to study the neural network models for rhythmic body movements during animal locomotion.
Such spinal neural networks are called central pattern generators (CPGs) that, 
as oscillator circuits, can generate coordinated signals for rhythmic body movements.
It has been well understood that {\it
``most neurons use action potentials (APs), brief and uniform pulses of electrical activity, to transmit information...
the strength at which an innervated muscle is flexed depends solely on the `firing rate,' the average number of APs per unit time (a `rate code')"} \cite{gerstner1997}. In other words, the transmitted signals are frequency modulated. 
For example, the frequency modulated  signals recorded in a real CPG network are shown in Fig.~\ref{fig:slug} 
for the dorsal-ventral swimmer Tritonia, a genus of sea slugs.

Frequency modulation (FM) is also one of most widely used technologies in engineering.
The conventional FM is the encoding of information in a carrier wave by varying the instantaneous frequency of the wave.
It has a long history with applications in  telecommunications, radio broadcasting, signal processing, and computing
\cite{saeedi2016synthetic,choi2017enhanced,wang2017parameter}.
For instance,  in radio transmission, FM, compared with  an equal power amplitude modulation  signal, 
has the advantage of a larger signal-to-noise ratio and therefore rejects radio frequency interference. 
Inspired by the aforementioned biological CPGs discovered in animal locomotion, engineers have applied 
artificial CPGs as oscillator circuits in robotics \cite{lewis2005cpg,ijspeert2007swimming,spaeth2020spiking}. 
Such artificial CPGs can generate frequency modulated signals,  resembling the firing bursts of APs,  to control robotic locomotion.

 However, it lacks mathematical principles, especially systems and control theories,  for frequency modulated multi-agent systems.
 A multi-agent system (MAS) refers to a group of multiple dynamic
entities (called agents) that share information or tasks to accomplish a common objective. 
The agents in an MAS can be circuits, autonomous vehicles, robots, power plants in electrical grid, 
and other dynamic entities. 
A CPG network in its dynamic formulation is a typical MAS
where each ganglion circuit (group of neurons) is regarded as an agent 
and they generate patterned signals through inter-agent coordination. 
In real biological CPGs, the signals of burst firing are frequency modulated
\cite{olivares2018,hachoumi2020,friesen:01}.
However, the mathematical principles for such CPGs have yet to be studied 
in terms of their coordinated behaviors. The research was mainly conducted based simplified 
unmodulated CPG models.
For example, the signals of burst firing 
are approximately represented by 
unmodulated sinusoidal-like waves when the adaptivity property of CPG circuits was studied
for Hirudo verbena, a species of leech. 
Such an unmodulated leech CPG model was used in \cite{chen:08,Iwasaki978} (and many other references) to
significantly reduce model complexity as there exists no systems and control theory
for frequency modulated MASs.

Over the past two decades, abundant researchers have
extensively studied  MASs and achieved numerous outcomes, 
across the disciplines of biology, statistical physics, computer sciences, sociology, management
sciences, and systems and control engineering.  
They have established many systems and control theories and engineering tools  
for various engineering tasks formulated as consensus, rendezvous, swarm, flocking, formation, etc. 
The early work on control of MASs focused on the so-called consensus problem where a networked set of agents 
intends to merge to a common state; see, e.g., \cite{jadbabaie:03,saber:04,ren:05} for simple linear homogenous MASs.
A major technical question was to understand the influence of network topology on collaborative behaviors.

The research has been widely extended to deal with more general homogeneous MASs in literature, e.g., 
 \cite{Scardovi2009,tuna2008lqr, Seo2009} for continuous-time systems and 
 \cite{tuna2008,youTAC2011,Movrica2013,stoorvogel2018}  for discrete-time systems.  
This line of research has also quickly spread to many other related problems such as synchronization, formation, flocking, swarming and rendezvous \cite{knorn2016,modares2019,jing2019}.  It is well known that
the collective property, in particular, consensus is closely related to the eigenvalue
distributions of the so-called Laplacian matrix, which in turn is crucially
influenced by the network topology.
 It is fair to say that the existing control strategies
for analysis and control of linear homogeneous MASs are mature using basic Laplacian matrix properties in graph theory.

The mature solutions to various control problems for homogenous MASs rely on 
the prerequisite that all the agents share share common dynamics. So, they mainly focus on trajectory synchronization. 
In the recent decade, researchers have put more efforts on synchronization of complicated 
 heterogeneous MASs, where the major challenge arises from the lack of common dynamics. 
To deal with this challenge, the work in \cite{Wieland2011,Kim2011} assumes
a common virtual exosystem that defines a
trajectory on which all agents synchronize to. 
Researchers also studied the problem using different concepts of homogenization in, e.g., 
\cite{yang2014output,Zhu2014}. Some other relevant techniques can be found in \cite{Lunze2012,Grip2013}, etc.

When researchers turned to nonlinear MASs, they first applied
the linear techniques to deal with nonlinearities under certain constraints on
the growth rates such as a globally Lipschitz-like condition; see, e.g., \cite{Yu2011,Su2011}.
The early research on nonlinear collective dynamic behaviors was also found on
synchronization of simply coupled nonlinear systems. The complexity is not from network but individual nonlinear
characteristics. 
Typical examples include synchronization of
master-slave chaotic systems \cite{Pecora1990}, Kuramoto models \cite{Chopra2009}, and  coupled harmonic oscillators \cite{Su2009}, etc. 

Synchronization of general nonlinear heterogeneous MASs is the state-of-the-art of this field. 
This problem was studied in \cite{Isidori2014,Chen:2014,Zhu2016} using a
generalized output regulation framework 
where the agents converge to various design of specified reference models. 
In a leader-following setting,  the problem was formulated as a cooperative
output regulation problem in, e.g., \cite{Huang2015} where agents achieve
synchronization on a trajectory that is assigned by one or more leaders. 
When the leader's system dynamics are unknown to some agents,
adaptive schemes were developed in, e.g.,  \cite{cai2017adaptive,JH_ADO,yan2019cooperative}. 
These design approaches for heterogeneous MASs can be illustrated as
``reference consensus + agent-reference regulation'' where the synchronized behavior is determined by 
the specified reference model dynamics (called virtual exosystem, exosystem, homogeneous model,
internal model, reference model, or leader, etc.) rather than inherent agent dynamics.

On the contrary, to adopt inherent agent dynamics rather than the specified 
common dynamics, people studied the so-called dynamics consensus/synchronization problem
as the prerequisite for trajectory synchronization. For example, 
the authors of  \cite{panteleydiss,panteley2017synchronization}
proposed the concept of emergent dynamics which is 
 an ``average'' of the ``units' drifts''.
An MAS reaches dynamic consensus if the agents' trajectories  converge to the one generated by the emergent dynamics.  
However, synchronization errors exist in these works due to 
the heterogeneity between unit dynamics  and emergent dynamics, 
which may be diminished by increasing the interconnection gain. 
The idea was also applied  on Andronov-Hopf oscillators \cite{maghenem2016singular}
and   Stuart-Landau oscillators  \cite{panteley2020practical}.
A more recent result on dynamics synchronization was developed in  \cite{yan2020,Hu2020}
in a new formulation of autonomous synchronization.
The result ensures that both the synchronized agent dynamics and the synchronized states are not specified a priori but autonomously determined by the inherent properties and the initial states
of agents.

In all the references discussed in the above survey of MASs, the final target of synchronization is on 
agent trajectories.  It can be achieved either through direct exchange of trajectory information via network (homogeneous MASs)
or indirectly through consensus of specified or autonomous reference models plus regulation of agent trajectories to references 
(heterogeneous MASs). A key stumbling block in generalizing the existing MAS control techniques  
to frequency modulated MASs is that the techniques require direct or indirect exchange and coordination of unmodulated values of agent trajectories in time domain.
However, for frequency modulated MASs, the target is not determined by  the instantaneous  values of agent trajectories.

The main contribution of this paper is the first attempt to establish a framework for synchronization of 
frequency modulated MASs.  
There is rare systems and control theory (analysis or controller synthesis) available for frequency modulated MASs.
The technical challenges lie in the fact that the frequency information is hidden from the
 instantaneous agent values  observed or transmitted in a network. 
 Analysis and control of  frequency modulated MASs
relies on an effective distributed frequency demodulation theory, which is still an open topic. 
Technically, we will first develop a new method to construct a frequency observer of a nonlinear system, especially subject to an external input representing synchronization error. 
None of the existing observer design methods is applicable in this scenario. 
Moreover, the gain from the external input to the observation error should be made arbitrarily small.
Secondly, we will design a frequency synchronization approach assuming the frequency observers work with 
appropriate observation errors.  As the two actions occur concurrently, we will also 
ensure a deliberately formulated small gain condition. 
The complete framework is established by an explicit controller architecture with rigorous mathematical proofs and 
demonstrated by numerical simulation.

The rest of the paper is organized as follows. The problem of synchronization 
of a frequency modulated MAS is formulated in Section~\ref{section:formulation}. 
A frequency observer for a nonlinear frequency modulated oscillator subject to network influence is designed in Section~\ref{section:observer}. Then, a frequency consensus algorithm is established 
 via  network  interaction subject to observation errors in Section~\ref{section:consensus}.
Based on the techniques for frequency observation and  frequency consensus
and examination of a small gain condition, the overall solution to the  synchronization problem 
is proposed in Section~\ref{section:synchronization}. The effectiveness of the 
approach is verified  in Section~\ref{section:example}. The paper is finally concluded in 
Section~\ref{section:conclusion}. 
 
\section{Problem Formulation} \label{section:formulation}

This paper is concerned about a network of  frequency modulated oscillators, or called agents, of the following dynamics
\begin{align}
\dot \sigma_i &= S \sigma_i +B\chi_i  \nonumber\\
\omega_i &= E\sigma_i +\omega_c  \nonumber\\
\dot x_i &= f(x_i)  \omega_i +f_o(x_i) ,\; i=1,\cdots, n \label{MAS}
\end{align}
where $\sigma_i\in \R^p$ is the state of the oscillator $i$. The information to be transmitted over a network
is $E\sigma_i$ modulated by a nonlinear carrier represented by the $x_i$-dynamics. In particular, 
the information $E\sigma_i$
is added on the top of the carrier's base frequency $\omega_c$ to form the frequency $\omega_i \in \R$ of the carrier.
Then, the actual transmitted signal over the network becomes $x_i \in \R^q$.
The matrices $S$, $B$, and $E$ are constant of appropriate dimensions. 
The vector functions $f$ and $f_o$ are continuously differentiable and specified a priori. 
In particular, the function $f$ satisfies  $f(x_i) \neq 0$ for any $x_i \neq 0$ so that $\omega_i$ is
effectively encoded in the signal $x_i$.

The objective is to design the input $\chi_i$ to each agent such that the states achieve synchronization 
in the sense of 
\begin{align} \label{target}
\lim_{t\rightarrow\infty} \| \sigma_i(t) -\sigma_j(t) \| =0,\; i, j = 1,\cdots, n,
\end{align}
or equivalently, in terms of the frequencies,
\begin{align} \label{target2}
\lim_{t\rightarrow\infty} \| \omega_i(t) -\omega_j(t) \| =0,\; i, j = 1,\cdots, n.
\end{align}
The design of $\chi_i$ relies on the received signals through a network.
Associated with the network,  we define a graph $\mathcal{G}=(\mathcal{V}, \mathcal{E}, A)$, where $\mathcal{V}=\left\{1, 2, \ldots, n \right\}$ indicates the vertex set, $\mathcal{E} \subseteq \mathcal{V} \times \mathcal{V}$ the edge set, and $A  \in \mathbb{R}^{n \times n}$ the associated adjacency matrix
whose  $(i,j)$-entry is $a_{ij}$. 
    For  $i, j \in \mathcal{V}$, $a_{i j} >0 $ if and only if $\left(i, j\right) \in \mathcal{E}$ and otherwise, $a_{i j}=0$. 
The set of neighbors of $i$ is denoted as $\mathcal{N}_{i}=\left\{j \; | \;\left(i, j\right) \in \mathcal{E}\right\}$.
  Let the $(i,j)$-entry of the Laplacian matrix  $L \in \mathbb{R}^{n \times n}$ be $l_{ij}(t)$ where $l_{ii}(t)=\sum_{j=1}^n a_{ij}(t)$ and $l_{ij}(t)=-a_{ij}(t)$  for $i\neq j$, $i, j \in \mathcal{V}$.

If the oscillator states $\sigma_i$, $i \in  \mathcal{V}$,  were available for transmission, 
synchronization can be easily achieved by a controller of the form
\begin{align} \label{chiideal}
\chi_i = -  M \sum_{j \in \N_i}  a_{ij} (\sigma_i -\sigma_j)  ,\; i \in  \mathcal{V} 
\end{align}
with a certain network connectivity condition for a properly designed matrix $M$.
However, in the present scenario, $\sigma_j$ is not available for agent $i$. The controller under 
consideration is modified to 
\begin{align}
\chi_i = & -  M \sum_{j \in \N_i}  a_{ij} (\sigma_i -\hat\sigma^i_j) \label{chii}\\
\hat\sigma^i_j =& \O(x_j),\; j \in \N_i,\; i \in  \mathcal{V},
\end{align}
where $\hat\sigma^i_j$ is the estimated value of $\sigma_j$ by agent $i$. The notation
$\O(x_j)$ is not necessarily a static function but  represents the output of a dynamic system whose input is $x_j$. 
The functionality of $\O(x_j)$ is to estimate the frequency $\omega_j$ and the oscillator 
state $\sigma_j$ using the received signal $x_j$. In this sense, $\O(x_j)$ is called a frequency observer. 
The specific design of $\O(x_j)$ for each agent is elaborated as follows, motivated by the dynamics \eqref{MAS},
\begin{align}
\dot {\hat \sigma}_j^i &= S \hat \sigma_j^i +\mu (x_j, \hat x_j^i, \hat\sigma_j^i) \nonumber\\
\hat\omega_j^i &= E\hat\sigma_j^i +\omega_c  \nonumber\\
\dot {\hat x}_j^i&=   f(  \hat x^i_j)\hat \omega_j^i +f_o(x_j)  +\kappa (x_j, \hat x_j^i ),\;
j\in \N_i , \; i \in  \mathcal{V} \label{observeri}
\end{align}
with $\hat x_j^i(0) = x_j(0)$, for some functions $\mu$ and $\kappa$.
Now, the aforementioned synchronization problem can be solved in a framework of three steps. 

Step 1 (Observer design):   Let $\tilde\sigma_j^i = \hat\sigma^i_j-\sigma_j$, $j\in \N_i$, $i \in  \mathcal{V}$, be the estimation errors. 
They can be put in a compact cumulative  form of  $\phi_i = \sum_{j \in \N_i} a_{ij}   \tilde\sigma^i_j$
and $\phi=\col (\phi_1,\cdots, \phi_n)$. Also, let $\chi=\col (\chi_1,\cdots, \chi_n)$. 
As the system \eqref{MAS} is non-autonomous with an external input $\chi_i$, 
the observer \eqref{observeri} is designed such that 
\begin{align*}
 \limsup_{t\rightarrow \infty}  \|\phi(t)  \| & \leq \gamma_\sigma   \limsup_{t\rightarrow \infty} \|\chi(t)\|
\end{align*}
for some gain $\gamma_\sigma$. 

Step 2 (Perturbed consensus): The control input $\chi_i$ to each agent is designed such that 
consensus is achieved with the estimation errors regarded as external perturbation in the sense of 
     \begin{align*} 
\limsup_{t\rightarrow \infty}  \|\chi (t) \|
\leq   \gamma_\chi \limsup_{t\rightarrow \infty}  \| \phi(t)\|
\end{align*} 
for some gain $\gamma_\chi$.

Step 3 (Small gain condition):   The two gains
 $\gamma_\sigma$ and $\gamma_\chi$ must satisfy a small gain condition to ensure the final target of 
 synchronization in terms of \eqref{target}, which puts an extra constraint on the design in the above two steps.
 
The framework of these three steps will be developed in the subsequent three sections, respectively. 
The technical challenges in the development are summarized as follows. 

(i) It requires a new method to construct a frequency observer of a nonlinear system, especially subject to an external input. 
None of the existing observer design methods is applicable in this scenario. 
Moreover, the gain from the external input to the observation error can be made arbitrarily small.

(ii) To decouple the complexity of the overall system, we first consider 
Step 1 and Step 2, separately. In particular,  Step 1 can be achieved under the assumption that $\chi(t)$ is bounded
and similarly, Step 2 can be achieved under the assumption that $\phi(t)$ is bounded.
It is a technical difficulty to ensure both assumptions can be satisfied concurrently 
in the closed-loop system through a certain interconnection. On satisfaction of these bounded assumptions, 
the proposed controller must also ensure a small gain condition, which requires deep understanding the coupling 
structure in the system dynamics. 
 
(iii) Among the agents, the communicated state is the modulated signal $x_i$, however, the target of the control strategy is 
synchronization of the state $\sigma_i$.  In fact, synchronization of the instantaneous values of 
 $x_i$ is not needed in this framework. 
So, the existing synchronization methods for nonlinear systems using a reference model (or an internal model) do not apply in the present setting.

\section{Design of a Frequency Observer} \label{section:observer}

 In this section, we study a frequency observer for a frequency modulated oscillator. 
 To simplify the notation, we omit the subscripts and superscripts in this section by focusing 
 on an individual scenario. More specially, we consider the model 
\begin{align}
\dot \sigma &= S \sigma  +B\chi   \nonumber\\
\omega  &= E\sigma +\omega_c  \nonumber\\
\dot x &=  f(x)  \omega +f_o(x)  \label{system}
\end{align}
and the corresponding observer
\begin{align}
\dot {\hat \sigma}&= S \hat \sigma +\mu(x, \hat x, \hat\sigma) \nonumber\\
\hat\omega &= E\hat\sigma+\omega_c  \nonumber\\
\dot {\hat x} &=f(\hat x) \hat \omega +f_o(x)   +\kappa (x , \hat x ) \label{observer}
\end{align}
with $\hat x(0) = x(0)$. 
 The following notations are used throughout the paper. 
For a symmetric real  matrix $P$, let  $\lambda_{\max}(P)$  and  $\lambda_{\min}(P)$ be the maximum and minimum  eigenvalues of $P$, respectively. And let $\varpi(P)=\frac{  \lambda_{\max}(P) }{ \lambda_{\min}(P)}$.

\blemma \label{lemma:sigma}
Consider the dynamics \eqref{system} and the observer \eqref{observer} with a constant $\omega_c >0$ and a time varying signal $\chi(t)$
satisfying $\|\chi(t)\| \leq b_\chi,\; \forall t\geq 0$ for a constant $b_\chi$. 
Suppose the frequency modulated signal $x(t)$ is bounded, i.e., $\underline \alpha \leq \| x(t) \| \leq \bar \alpha, \forall t\geq 0$
for some constants  $\underline\alpha, \bar\alpha>0$. 
For any constants $b_o, \gamma >0$,  
if there exist matrices $K_o$ and  $P=P\t >0$ such that \begin{align}
P (S  -K_o E ) +(S  -K_o E )\t P & =-I  \nonumber\\ 
\|PB\| \sqrt{\varpi(P)}&<\gamma_o \label{PSSP}
 \end{align} 
 for $\gamma_o=\gamma/8$,
 then there exist control functions $\mu$ and $\kappa$ such that
\begin{align}
\limsup_{t\rightarrow \infty} \| \hat \sigma(t) -\sigma(t)\| & \leq \gamma \limsup_{t\rightarrow \infty} \|\chi(t)\|  \label{sigmachi}
\end{align}
for  \begin{align}
\| \hat \sigma(0) -\sigma(0)\| & \leq b_o. \label{boundsigmax}
\end{align}
\elemma

\bproof For the convenience of presentation, denote the observation errors
$\tilde\sigma =\hat \sigma -\sigma$ and $\tilde x =\hat x - x$. 
With the relationship
$\hat \omega  - \omega  =
E\tilde \sigma$, the error dynamics can be written as 
\begin{align*}
\dot {\tilde \sigma} &= S \tilde \sigma  +\mu (x, \hat x, \hat\sigma) -B\chi \\
 \dot {\tilde x} &=f(\hat x) \hat \omega  -f(x)  \omega    +\kappa (x , \hat x ) \\
 &=  f(x) E\tilde \sigma   +[f(\hat x) -f(x)] \hat \omega +\kappa (x, \hat x ).
\end{align*}
 
 Pick  a matrix function 
 \begin{align}
 K(  \hat x)= K_o \frac{ f\t( \hat  x) }{\| f( \hat  x)  \|^2 }.  \label{Kdef}
   \end{align}
 It will be proved later that $\| f( \hat  x(t))  \| >0, \;\forall t\geq 0$, to validate the definition of $K$.
Direct calculation shows that 
 \begin{align*}
   K(\hat x) f( \hat x)  
&=
    K_o \frac{ f\t(\hat x) }{\| f(\hat x)  \|^2 } f(\hat x) = K_o . \end{align*}
Then, one has  
\begin{align*}
S  \tilde \sigma 
 - K(\hat x) f( \hat x)E   \tilde \sigma 
 = S  \tilde \sigma -K_o E   \tilde \sigma 
= 
 \bar S  \tilde \sigma 
 \end{align*}
 where $\bar S=S  -K_o E$ is Hurwitz with \eqref{PSSP} that can be rewritten as 
  \begin{align}
P \bar S +\bar S\t P & =-I \nonumber \\
8\|PB\| \sqrt{ \varpi(P)} & <\gamma.\label{PBlambda}
 \end{align}

We introduce a new variable 
\begin{align} \label{nu}
\nu = {\tilde \sigma} - K(\hat x)  {\tilde x},
\end{align}
that gives
\begin{align*}
\dot \nu = & S \tilde \sigma  +\mu (x, \hat x, \hat\sigma) -B\chi  \\
& - K(\hat  x) [  f(x) E\tilde \sigma  +[f(\hat x) -f(x)] \hat \omega  +\kappa (x, \hat x ) ] \\
&- \dot K (\hat  x)  {\tilde x}     \\
 = & S \tilde \sigma 
 - K( \hat  x)f(x) E\tilde \sigma  +\bar\mu (x, \hat x )   -B\chi  \\
  = & S \tilde \sigma 
  - K( \hat  x)f(\hat x) E\tilde \sigma
  - K( \hat  x)f(x) E\tilde \sigma \\
& +  K( \hat  x)f(\hat x) E\tilde \sigma
   +\bar\mu (x, \hat x )   -B\chi  \\
  = &
  \bar S \tilde \sigma 
   +  K( \hat  x) [f(\hat x) -f(x)] E\tilde \sigma
   +\bar\mu (x, \hat x )   -B\chi  
    \end{align*}
where
\begin{align*}
\mu (x, \hat x, \hat\sigma) 
=&   K( \hat  x)\kappa (x, \hat x)  +K( \hat  x) [f(\hat x) -f(x)] \hat \omega \\ & +K'(x, \hat x, \hat\sigma)  {\tilde x}  
  +\bar\mu (x, \hat x ) \\
\bar\mu (x, \hat x)   = & -\bar S   K(\hat x)  {\tilde x} 
- K( \hat  x) [f(\hat x) -f(x)] E  K(\hat x)  {\tilde x} 
 \end{align*}
 and      \begin{align*}
K'(x, \hat x, \hat\sigma)= & K_o \left(  \frac{ \| f( \hat  x)  \|^2  I_q
- 2 f(\hat x)  f\t(\hat x)   
 } {\| f( \hat  x)  \|^4 }  \times \right. \\
 & \left. \frac{\partial f( \hat  x) }{ \partial \hat x} 
 \left(f(\hat x) \hat \omega +f_o(x)   +\kappa (x , \hat x ) \right)
 \right)\t .
   \end{align*}
 It is noted that, along the trajectory of \eqref{observer},      \begin{align*}
\dot K(\hat x) =\frac{\partial   K(\hat x)}{\partial \hat x}
\dot {\hat x}
=K'(x, \hat x, \hat\sigma) .
   \end{align*}
The $\nu$-dynamics can be further simplified as follows
\begin{align}
\dot \nu      = & \bar S  \nu   +  K( \hat  x) [f(\hat x) -f(x)] E \nu     -B\chi . \label{nudynamics}
\end{align}

Next, noting $E\tilde \sigma =E \nu +E K(\hat x)  {\tilde x}$, 
the $\tilde x$-dynamics can be rewritten as 
\begin{align}
 \dot {\tilde x} &=  f(x) E\nu -  \beta \tilde x  \label{tildexdynamics}
\end{align}
where 
\begin{align*}
\kappa (x, \hat x )  = -\beta \tilde x  -  f(x) E K(\hat x)  {\tilde x}- [f(\hat x) -f(x)] \hat \omega\end{align*}
and $\beta$ is to be specified later.

As $f$ is continuously differentiable, 
there exist positive numbers $\alpha,\theta>0$ such that (see Lemma~11.1 of \cite{chen2015stabilization})
  \begin{align*}
  \| [f(x +\tilde x) -f(x)] E \nu\|  & \leq \theta \|\nu \| \|\tilde x\|\\
   \|f(x) E\nu \|& \leq \theta  \alpha \|\nu\|,\; \forall \underline\alpha \leq \| x \| \leq \bar \alpha,\;  \|\tilde x\| \leq \underline\alpha/2.
    \end{align*}
Also, we define positive numbers $b_x, b_K >0$ as
  \begin{align} \label{bx}
  b_x  = \arg\max_{b\leq \underline\alpha/2}\{8\theta \| P K(x+\tilde x) \| \|  \tilde x \| \leq  1,\; \forall  \underline\alpha \leq \| x \| \leq \bar \alpha, \nonumber \\ \|\tilde x\| \leq b \}
 \end{align}
 and 
    \begin{align}
  b_K=  \max_{\underline\alpha \leq \| x \| \leq \bar \alpha, \|\tilde x\| \leq b_x} \| K(x+\tilde x)  \|. \label{bK}
 \end{align}
From the above definition, one has $b_x \leq \underline \alpha/2$.
For $\underline\alpha \leq \| x \| \leq \bar \alpha, \|\tilde x\| \leq b_x$, one has $\underline\alpha/2 \leq \|\tilde x +x\| \leq \bar \alpha + \underline\alpha/2$, 
so $b_K$ is finite according to the definition of the function $K$. Then, let $\beta$ satisfy the following conditions
   \begin{align}\label{beta}
      \beta  &>  \frac{\varpi(P) b_o^2+ 8\rho^2 b_\chi^2 \varpi(P)}{4b_x^2} +\theta^2 \alpha^2 \\
      \beta & >  \frac{b_K^2}{4}+\theta^2 \alpha^2  \label{beta2}\\
       \beta & \geq \frac{1}{4\lambda_{\min}(P)}+ \lambda_{\min}(P)\theta^2 \alpha^2. \label{beta3}
  \end{align}
    The remaining analysis will be on the interconnected system composed of \eqref{nudynamics} and \eqref{tildexdynamics}.

 Next, we will prove
    \begin{align}
   \|\tilde x(t)\| < b_x, \;\forall t \geq 0. \label{xbound}
     \end{align}
If \eqref{xbound} does not hold, there exists a finite $T >0$ such that
$\|\tilde x(t)\| < b_x, \;\forall t\in [0, T)$ and $\|\tilde x(T)\| = b_x$, 
because $\|\tilde x(0)\| =0$ and $\tilde x(t)$ is a continuous function of $t$.
A contradiction will be derived below.

 First, for the system \eqref{nudynamics}, pick a Lyapunov function candidate
 $V_1(\nu) = \nu\t P \nu/2$, whose time derivative satisfies 
  \begin{align*}
\dot V_1(\nu) = & - \| \nu\|^2  +2  \nu\t  P K(x+\tilde x ) [f(\hat x) -f(x)] E \nu      -2 \nu\t  P B\chi  \\
\leq & - \| \nu\|^2  +2 \theta  \| P K(x+\tilde x ) \| \|\tilde x\|    \| \nu\|^2  \\ & +\frac{1}{4} \|\nu\|^2 +  4\|PB\|^2 \|\chi\|^2.
 \end{align*}
For  $t\in [0, T]$, one has $\|\tilde x(t)\| \leq b_x$ and hence
$8\theta \| P K(x+\tilde x) \| \|  \tilde x \| \leq  1$. Therefore, 
 \begin{align}\label{dV10}
\dot V_1(\nu) 
\leq & - \| \nu\|^2  +\frac{1}{4}  \| \nu\|^2   +\frac{1}{4} \|\nu\|^2 +  4\|PB\|^2 \|\chi\|^2  \nonumber\\
\leq & - \frac{1}{2}  \| \nu\|^2  +  4\|PB\|^2 \|\chi\|^2 
 \end{align}
and 
 \begin{align}
\dot V_1(\nu) \leq   - \frac{1}{ \lambda_{\max}(P)}   V_1(\nu)  +  4\|PB\|^2 \|\chi\|^2. \label{dV1}
 \end{align}
From \eqref{dV1}, the comparison principle gives
 \begin{align*}
 & V_1(\nu(t))\\
  \leq &   V_1(\nu(0)) + 4\|PB\|^2  \int_0^t \exp (- \frac{1}{ \lambda_{\max}(P)} (t-\tau) )
  \|\chi(\tau)\|^2 d\tau  \\
  \leq  & V_1(\nu(0)) + 4\|PB\|^2 b_\chi^2 \int_0^t \exp (- \frac{1}{ \lambda_{\max}(P)} (t-\tau)) d\tau\\
= &    V_1(\nu(0)) + 4\|PB\|^2 b_\chi^2 \lambda_{\max}(P) \left[1- \exp (- \frac{1}{ \lambda_{\max}(P)} t)\right]\\
 \leq& V_1(\nu(0)) + 4\|PB\|^2 b_\chi^2 \lambda_{\max}(P) \\
  \leq & \lambda_{\max}(P) \|\nu(0)\|^2/2 + 4\|PB\|^2 b_\chi^2 \lambda_{\max}(P).   \end{align*}
 From the above calculation, one has
 \begin{align*}
 \lambda_{\min}(P)\|\nu(t)\|^2/2  &\leq V_1(\nu(t)) \\ &  \leq   
 \lambda_{\max}(P) \|\nu(0)\|^2/2+ 4\|PB\|^2 b_\chi^2 \lambda_{\max}(P),  \end{align*}
and then
\begin{align} \label{nubound}
 \|\nu(t)\|^2 \leq \varpi(P) \|\nu(0)\|^2+ 8\|PB\|^2 b_\chi^2 \varpi(P).  \end{align}
By \eqref{boundsigmax} and the fact
\begin{align*}
\nu(0) = {\tilde \sigma}(0) - K(\hat x(0))  {\tilde x(0)} = {\tilde \sigma}(0),
\end{align*}
\eqref{nubound} implies 
 \begin{align}\label{nubound2}
 \|\nu(t)\|^2 \leq \varpi(P) b_o^2+ 8\|PB\|^2 b_\chi^2 \varpi(P). \end{align}

Secondly, for the system  \eqref{tildexdynamics}, pick a Lyapunov function candidate
 $V_2(\tilde x) = \|\tilde x\|^2 /2$, whose time derivative satisfies 
  \begin{align}\label{dV2}
\dot V_2(\tilde x)= & -\beta \|\tilde x\|^2+ \tilde x\t f(x) E \nu  \nonumber  \\
\leq& -\beta \|\tilde x\|^2 + \theta \alpha \|\tilde x\| \|\nu\| \nonumber \\
\leq& -(\beta - \theta^2 \alpha^2  ) \|\tilde x\|^2  + \frac{1}{4}\|\nu\|^2\nonumber \\
 \leq& -2(\beta - \theta^2 \alpha^2) V_2(\tilde x)  + \frac{1}{4}\|\nu\|^2.
 \end{align}
 It is noted that  $\beta - \theta^2 \alpha^2 >0$  from \eqref{beta} or \eqref{beta2}.
Applying the comparison principle again gives
  \begin{align*}
&  V_2(\tilde x(t)) \\ \leq&     V_2(\tilde x(0)) + \frac{1}{4}  \int_0^t \exp (-2(\beta - \theta^2\alpha^2  ) (t-\tau) )
  \|\nu(\tau)\|^2 d\tau  \\
   \leq&      \frac{\varpi(P) b_o^2+ 8\|PB\|^2 b_\chi^2 \varpi(P)}{4} \int_0^t \exp (-2(\beta - \theta^2\alpha^2 ) (t-\tau) )
 d\tau  \\
  \leq&      \frac{\varpi(P) b_o^2+ 8\|PB\|^2 b_\chi^2 \varpi(P)}{8(\beta - \theta^2\alpha^2 )}  
  [1-\exp (-2(\beta - \theta^2\alpha^2 )t ]\\
\leq&      \frac{\varpi(P) b_o^2+ 8\|PB\|^2 b_\chi^2 \varpi(P)}{8(\beta - \theta^2 \alpha^2 )},
    \end{align*}
 where \eqref{nubound2} is used in the above calculation. By \eqref{beta}, one has
   \begin{align*}
 \|\tilde x(t)\|^2 =2 V_2(\tilde x(t)) 
\leq&      \frac{\varpi(P) b_o^2+ 8\rho^2 b_\chi^2 \varpi(P)}{4(\beta - \theta^2 \alpha^2)} <   b_x^2    \end{align*}
for $t\in [0, T]$. It is a contradiction with $\|\tilde x(T)\| = b_x$.
From the above argument, one has \eqref{xbound} proved. Also, 
 $\|K(\hat x(t))\|   \leq b_K$,  $\hat x(t) \neq 0$,  $f( \hat  x(t)) \neq 0$, 
 \eqref{dV10}, and \eqref{dV2} hold for all $t\geq 0$.
The definition of $K$ in \eqref{Kdef} is thus validated.

Next, as \eqref{beta2} is equivalent to 
     \begin{align*} \sqrt{\frac{1}{4(\beta - \theta^2 \alpha^2) }} <\frac{1}{b_K  },
     \end{align*}  \eqref{dV2}  implies that
the $\tilde x$-system \eqref{tildexdynamics} is input-to-state stable viewing $\tilde x$ as the state
and $\nu$ as the input and particularly (see Theorem~2.7 of \cite{chen2015stabilization}),
     \begin{align*}
 \limsup_{t\rightarrow \infty} \|  {\tilde x(t)}\|  \leq\frac{1}{b_K  } \limsup_{t\rightarrow \infty} \|\nu(t)\|.
    \end{align*}
And the following inequalities hold
  \begin{align}\limsup_{t\rightarrow \infty}
  \|K(\hat x(t))  {\tilde x(t)}\| \leq b_K\limsup_{t\rightarrow \infty} \|  {\tilde x(t)}\|  \leq \limsup_{t\rightarrow \infty} \|\nu(t)\|.
  \label{limK}
   \end{align}

For $\bar V_2(\tilde x) =\lambda_{\min}(P) \|\tilde x\|^2 /2$, a variant of  \eqref{dV2}  is written as follows
  \begin{align}
\dot {\bar V}_2(\tilde x)= & -\lambda_{\min}(P)\beta  \|\tilde x\|^2+\lambda_{\min}(P) \tilde x\t f(x) E \nu  \nonumber \\
\leq& -\lambda_{\min}(P)\beta \|\tilde x\|^2 +\lambda_{\min}(P) \theta \alpha \|\tilde x\| \|\nu\| \nonumber\\
\leq& -(\lambda_{\min}(P)\beta - \lambda^2_{\min}(P)\theta^2 \alpha^2) \|\tilde x\|^2  + \frac{1}{4}\|\nu\|^2.
\label{dbarV2}
 \end{align}
 It is ready to construct the composite Lyapunov function 
  \begin{align*}
  V(\nu,\tilde x) =  V_1(\nu)  + \bar V_2(\tilde x) 
\end{align*}
that satisfies
 \begin{align*}
 \lambda_{\min}(P)\|\col(\nu, \tilde x)\|^2 /2 \leq & V(\nu,\tilde x)\\ \leq & \lambda_{\max}(P) \|\nu\|^2/2+ 
  \lambda_{\min}(P)\|\tilde x\|^2 /2 \\ \leq &  \lambda_{\max}(P)\|\col(\nu, \tilde x)\|^2/2.
  \end{align*}
Using  \eqref{dV10} and \eqref{dbarV2}, the derivative of $V(\nu,\tilde x)$ satisfies
 \begin{align*}
\dot V(\nu,\tilde x)
\leq & - \frac{1}{4}  \| \nu\|^2   -(\lambda_{\min}(P)\beta - \lambda^2_{\min}(P)\theta^2 \alpha^2) \|\tilde x\|^2   \\
&+  4\|PB\|^2 \|\chi\|^2  \\
\leq & - \frac{1}{4} \|\col(\nu, \tilde x)\|^2    +  4\|PB\|^2 \|\chi\|^2 
 \end{align*}
 noting 
$\lambda_{\min}(P)\beta - \lambda^2_{\min}(P)\theta^2 \alpha^2\geq 1/4$
according to \eqref{beta3}.
As a result, 
the $(\nu,\tilde x)$-system composed of \eqref{nudynamics} and  \eqref{tildexdynamics} is input-to-state stable viewing $(\nu,\tilde x)$ as the state
and $\chi$ as the input. In particular, one has
 \begin{align}
\limsup_{t\rightarrow \infty} \| \col(\nu(t), \tilde x(t)) \| & \leq \gamma/2 \limsup_{t\rightarrow \infty} \|\chi(t)\|   \label{sigmachi2}
\end{align}
as \eqref{PBlambda} is equivalent to 
 \begin{align*}
\frac{\gamma}{2} >  \sqrt {\frac{\lambda_{\max}(P)4\|PB\| ^2  }{ \lambda_{\min}(P) /4 } } =4\|PB\| \sqrt{\varpi(P)}.
  \end{align*}
Finally, using \eqref{nu}, \eqref{limK}, and \eqref{sigmachi2}, the following 
conclusion is verified
 \begin{align*}
& \limsup_{t\rightarrow \infty} \|\tilde \sigma(t)\| \\ \leq &
\limsup_{t\rightarrow \infty}\|\nu(t)\| +\limsup_{t\rightarrow \infty} \|K(\hat x(t))  {\tilde x(t)}\|\\
\leq  &
2\limsup_{t\rightarrow \infty}\|\nu(t)\| \leq 2 \limsup_{t\rightarrow \infty} \| \col(\nu(t), \tilde x(t)) \| \\ \leq &
 \gamma \limsup_{t\rightarrow \infty} \|\chi(t)\|,
\end{align*}
which is \eqref{sigmachi}. 
The proof is thus completed. 
\eproof

\bremark Design of a state observer for estimating  the internal state of a given system is an
important topic in control theory. For example, there are many practical applications of a Luenberger observer and a Kalman filter.
Researchers have also studied many nonlinear observer design methods. 
The observer in \cite{xia1989nonlinear} is based on the observer error linearization approach. 
An extended Luenberger observer was designed in \cite{ZEITZ1987149} using the technique of linearization and 
transformation into a nonlinear observer canonical form. 
High gain observers were studied in many references, e.g.,  \cite{khalil2014high,Astolfi2015}, 
for nonlinear systems in canonical observability form.  Neither the approach based on linearization nor a
canonical form does not apply in the present paper  due to the special structure of the system under consideration. 
However, the new approach takes advantage of the interconnection structure between the $\sigma$ and $x$ subsystems
and derives a Lyapunov function based argument. 
In addition, an arbitrarily specified gain from an external input to the estimation error further complicates the problem. 
\eremark

\bremark 
The gain from $\|\chi\|$ to the estimation error $\|\hat \sigma -\sigma\|$ is characterized in the lemma by an
appropriately designed observer. It is noted that the influence of $\chi$ could be completely removed by compensation in 
the observer  \eqref{observer} such as   $\dot {\hat \sigma}= S \hat \sigma +B\chi +\mu(x, \hat x, \hat\sigma)$. 
However, such compensation is not  applicable when the observer is used for network synchronization in this paper. 
Consider the observer \eqref{observeri} of agent $i$ to estimate the frequency of agent $j$. The only signal agent $i$ is able to access is the state $x_j$, not the input $\chi_j$ of agent $j$. 
\eremark

\bremark It is assumed that the frequency modulated signal $x(t)$ is bounded, i.e., $\underline \alpha \leq \| x(t) \| \leq \bar \alpha, \forall t\geq 0$,
in the lemma. It is a reasonable assumption for a persistently exciting oscillator. For example, for an oscillator
 \begin{align*}
\dot x = \left[\begin{array}{cc} 0 & \omega(t) \\-\omega(t) & 0 \end{array} \right] x,
 \end{align*}
one has $\| x(t) \| =\|x(0)\|,\; \forall t\geq 0$ that is bounded no matter how $\omega(t)$ varies, because 
 \begin{align*}
\frac{d \| x(t)\|^2}{dt} =
2 x\t \dot x =2 x\t  \left[\begin{array}{cc} 0 & \omega(t) \\-\omega(t) & 0 \end{array} \right] x=0. 
 \end{align*}
Oscillation is also typically generated as a stable limit cycle which is a closed trajectory and bounded. 
\eremark

\medskip

To close this section, we discuss the existence of the solution to \eqref{PSSP}
for a special case. 

\bcorollary \label{corK}
For 
\begin{align*}
  S=\left[\begin{array}{cc} 0 & \varsigma \\-\varsigma & 0 \end{array} \right],\; \varsigma>0,\;
E=\left[\begin{array}{cc}  1 &  0 \end{array} \right],\;
 B=\left[\begin{array}{c}  1 \\ 1 \end{array} \right],  
 \end{align*}
and  any constant $\gamma_o>0$, there exist matrices $K_o$ and  $P=P\t >0$ such that
\eqref{PSSP} is satisfied. 
    \ecorollary

\bproof For a positive $\rho$, let  $K_o= \rho B /EB +SB /EB$.
Then, 
    \begin{align*}
\bar S   =& S  -K_o E 
 = S - \rho BE /EB - SBE /EB \\
=& S (I-BE /EB) + \rho B E/EB\\
=&\left[\begin{array}{cc} -\varsigma-\rho  & \varsigma \\ -\rho & 0 \end{array} \right].
  \end{align*}
 We can solve the equation in \eqref{PSSP} as
       \begin{align*}
P  =\frac{1}{2\varsigma}  \left[\begin{array}{cc}   1 & -1  \\ -1 & 
\frac{\rho+2\varsigma}{ \rho}   \end{array} \right]
  \end{align*}
which has two eigenvalues $-\varsigma$ and $-\rho$ and hence is Hurwitz. 
 The characteristic equation of $P$ is
          \begin{align*}\det (P_o-\lambda I) 
 =&\lambda^2 -  \frac{\rho+\varsigma}{ \rho \varsigma}  \lambda 
 + \frac{1}{2 \rho \varsigma}  =0,
  \end{align*}
  which implies  
       \begin{align*}
 \varpi(P)+\varpi^{-1}(P) 
 = \frac{  \lambda_{\max}(P) }{ \lambda_{\min}(P)}+\frac{  \lambda_{\min}(P) }{ \lambda_{\max}(P)}
  =  \frac{2\rho+2\varsigma^2/\rho+2 \varsigma}{  \varsigma}
    \end{align*}
  and hence
        \begin{align*}
 \varpi(P)< \frac{2\rho+2\varsigma^2/\rho+2 \varsigma}{  \varsigma}. 
    \end{align*}
  It is noted that $|PB\| =1/\rho$.
Then,         \begin{align*}
\|PB\|^2  \varpi(P) &< \frac{2 +2\varsigma^2/\rho^2+2 \varsigma/\rho}{ \rho \varsigma}
\leq \frac{2 +2\varsigma^2+2 \varsigma}{ \rho \varsigma} 
\leq \gamma_o^2
    \end{align*}
for 
   \begin{align*}
\rho\geq \max\{1, \frac{2 +2\varsigma^2+2 \varsigma}{\gamma_o^2 \varsigma}   \}.
    \end{align*}
The proof is thus completed.   
 \eproof

\section{Perturbed Consensus} \label{section:consensus}

Consider the $\sigma_i$-dynamics of the system \eqref{MAS}, that is,
\begin{align}
\dot \sigma_i &= S \sigma_i +B\chi_i   ,\; i \in  \mathcal{V}.  \label{dyn-sigma}
\end{align}
If the state $\sigma_i$ could be directly transferred over the network, without 
frequency modulation, the synchronization target \eqref{target} can be easily achieved 
through the ideal controller \eqref{chiideal}.  However, with frequency modulation, the controller takes 
the form \eqref{chii} where the estimated state $\hat \sigma_j^i$ rather than $\sigma_j$ is available 
for the design.  
The controller \eqref{chii} can be rewritten as the ideal controller \eqref{chiideal} subject to perturbation (i.e., 
the estimation error),
\begin{align}
\chi_i  
= -M \sum_{j \in \N_i}  a_{ij} (\sigma_i - \sigma_j)  + M \sum_{j \in \N_i}  a_{ij}  \tilde\sigma^i_j, \;
i \in  \mathcal{V}. \label{chi_b}
\end{align}
 For $\phi_i = \sum_{j \in \N_i} a_{ij}   \tilde\sigma^i_j$, the system composed of \eqref{dyn-sigma}
and \eqref{chi_b} becomes
\begin{align} \label{dsigmai}
\dot \sigma_i = S \sigma_i  -BM \sum_{j \in \N_i}  a_{ij} (\sigma_i - \sigma_j)  +BM \phi_i , \; i \in  \mathcal{V}. 
\end{align}
For the convenience of presentation, let 
$\phi =\col(\phi_1,\cdots, \phi_n)$ and
$\sigma =\col(\sigma_1,\cdots, \sigma_n)$.

Let $\br, \bone\in\mathbb{R}^{n}$ be the left and right eigenvectors corresponding to the
eigenvalue zero of the  Laplacian matrix  $L$
and satisfying $\br\t\bone=1$. There exist matrices $W\in\mathbb{R}^{(n-1)\times n}$
and $U\in\mathbb{R}^{n\times(n-1)}$ such that
\begin{equation}
T=\left[\begin{array}{c}
\br\t\\
W
\end{array}\right],\;  T^{-1}=\left[\begin{array}{cc}
\bone & U\end{array}\right],\label{eq:T}
\end{equation}
 $W\boldsymbol{1}=0$ and $\br\t U=0$. Then
\begin{equation}
T L  T^{-1}=\left[\begin{array}{cc}
0 & 0\\
0 & H 
\end{array}\right]\label{eq:H}
\end{equation}
where 
\begin{align}
H  =W L  U
\end{align}
is called an H-matrix associated with $L$. It is well known that, under Assumption~\ref{ass-connect} given below, all the eigenvalues 
of $H$ have positive real parts. 

\basm \label{ass-connect}
The directed graph $\mathcal{G}$ contains a spanning tree.
\easm

Now, the main result about perturbed consensus is given in the following lemma. 

\blemma \label{lemma:chi}
Consider the system \eqref{dyn-sigma} with \eqref{chi_b}
under Assumption~\ref{ass-connect}, if $M$ is selected such that
 \begin{align}
 A_{\zeta} = I_{n-1} \otimes    S - H \otimes BM \label{M}
 \end{align}
  is Hurwitz, then 
\begin{align} \label{phi2chi}
\limsup_{t\rightarrow \infty}  \|\chi (t) \|
\leq   \gamma_\chi \limsup_{t\rightarrow \infty}  \| \phi(t)\|
\end{align} 
for some constant $\gamma_\chi$ and any bounded $\phi(t)$.

\elemma

\bproof  The the system composed of \eqref{dyn-sigma}
and \eqref{chi_b}, i.e.,  \eqref{dsigmai}, can be put in a compact form
\begin{align}
\dot \sigma =& (I_n \otimes S -  L \otimes BM) \sigma   +(I_n \otimes BM) \phi. \label{dsigma}
\end{align}
We introduce the following coordinate transformation
\begin{align}
\left[\begin{array}{c}
\bar \sigma \\
\zeta
\end{array}\right]=\left(T \otimes I_p \right) \sigma, \;
\sigma=\left(T^{-1} \otimes I_p \right) \left[\begin{array}{c}
\bar \sigma\\
\zeta
\end{array}\right]. \label{eq:x_transform}
\end{align}
 In the new coordinate,  the system \eqref{dsigma} is equivalent to
 \begin{align*}
\left[\begin{array}{c}
\dot {\bar \sigma} \\
\dot \zeta
\end{array}\right]  =&   (I_n \otimes S -T   L T^{-1}  \otimes BM)   \left[\begin{array}{c}
\bar \sigma \\
\zeta
\end{array}\right]   +(T  \otimes BM) \phi,
\end{align*}
that is
\begin{align}
\dot{\bar \sigma} & =  S \bar{\sigma} +(\br\t  \otimes BM) \phi \nonumber\\
\dot{\zeta} & =  A_{\zeta}\zeta  +(W  \otimes BM) \phi. \label{zeta}
\end{align}

 As $A_{\zeta}$  is Hurwitz, there exists $Q = Q\t >0$ as the solution to the Lyapunov equation
\begin{align} Q A_{\zeta} +A_{\zeta}\t Q =-2 I. \label{lya.equ}
\end{align}
The derivative of the function
\begin{align} V(\zeta )= \zeta\t   Q \zeta  \label{Vzeta}
\end{align}
along the trajectories of \eqref{zeta} is
\begin{align} \dot V(\zeta ) = & - 2 \|\zeta\|^2 + 2 \zeta\t   Q (W  \otimes BM) \phi \nonumber\\
= & - \|\zeta\|^2 + \| Q (W  \otimes BM)\|^2  \| \phi\|^2. \label{dVzeta}
\end{align}
It implies that the $\zeta$-system \eqref{zeta} is input-to-state stable viewing $\zeta$ as the state
and $\phi$ as the input. In particular, 
\begin{align}
\limsup_{t\rightarrow \infty} \| \zeta(t)\| & \leq \gamma_\zeta \limsup_{t\rightarrow \infty} \|\phi(t)\|    \label{limzeta}
\end{align}
for
\begin{align} \label{gammazeta}
\gamma_\zeta =\sqrt{\varpi(Q )} \| Q (W  \otimes BM)\|.
\end{align}

 Next,  from the relationship
\begin{align*}
  ( L \otimes M ) \sigma = ( L T^{-1} \otimes M ) \left[\begin{array}{c}
\bar \sigma\\
\zeta
\end{array}\right] =( L U \otimes M )\zeta
\end{align*}
one has
\begin{align*}
\limsup_{t\rightarrow \infty} \|  ( L \otimes M )  \sigma(t)\| & \leq
\|  L  U \otimes M  \|  \limsup_{t\rightarrow \infty} \|\zeta(t)\|\\ &   \leq
\| L U \|\|M\|\gamma_\zeta \limsup_{t\rightarrow \infty} \|\phi(t)\|.
\end{align*}

Finally,  \eqref{chi_b} can be put in the following compact form
\begin{align}
\chi  
= -( L \otimes M)  \sigma    + (I_n \otimes M) \phi
\end{align}
that gives
\begin{align*}
& \limsup_{t\rightarrow \infty}  \|\chi (t) \| \\
\leq &  \| L U \|\|M\|\gamma_\zeta \limsup_{t\rightarrow \infty} \|\phi(t)\|  + \|I_n \otimes M\| \limsup_{t\rightarrow \infty}  \| \phi(t)\| \\
= & \|M\| (\| L U \| \gamma_\zeta +1)  \limsup_{t\rightarrow \infty}  \| \phi(t)\|.
\end{align*}
The proof is thus complete for $\gamma_\chi  = \|M\| (\| L U \| \gamma_\zeta +1)$. \eproof

\bremark
Under Assumption~\ref{ass-connect},  let
 $\lambda_i,\;i=2,\cdots n$, be the eigenvalues of $H$ that have
 positive real parts.  If $(S,B)$ is stabilizable, 
 let $G=G\t>0$ be the solution to the following algebraic Riccati equation
 \begin{align} S\t G+GS-\lambda^{*}GB B\t G+\epsilon I=0 
  \end{align}
for some $0 < \lambda^{*} \leq \min_{i=2,\cdots N}\Re\{\lambda_{i}\}$ and
 $\epsilon>0$. Let the controller gain matrix be $M=B\t G/2$.
Then, all the eigenvalues of $S-\lambda_i B M$ have strictly negative real parts and hence
$A_\zeta$ in \eqref{M} is Hurwitz. \eremark

\bremark
The concept of perturbed consensus was originally proposed in \cite{Zhu2016}
to deal with synchronization of  heterogeneous nonlinear MASs using output communication, 
where the perturbation represents some trajectory regulation error. In this reference, perturbed consensus is 
formulated as a constant gain from the perturbation $\phi$ to the consensus error $\|\sigma- \bone\otimes \sigma_o\|$
for some agreed trajectory $\sigma_o$.
In this paper, to facilitate the subsequent analysis of the overall system, perturbed consensus is formulated 
in terms of a gain from the perturbation $\phi$ to $\chi$, i.e., \eqref{phi2chi}. In Lemma~\ref{lemma:chi},
perturbed consensus is achieved under the assumption 
that  the perturbation $\phi(t)$ is always bounded for $t \geq 0$, which will be ensured when the overall system is considered
in the following section. 
 \eremark

\section{Synchronization of a Frequency Modulated MAS} \label{section:synchronization}

With the frequency observers and the controllers for perturbed consensus proposed in the previous two sections, 
it is ready to obtain the result of synchronization of a frequency modulated MAS using a small gain based argument. The result is stated in the following theorem.

\btheorem \label{thm-FMMAS}
Consider the dynamics \eqref{MAS}, the controller \eqref{chii}, and the observer \eqref{observeri} with a constant $\omega_c >0$
under Assumption~\ref{ass-connect}.
Suppose the frequency modulated signal $x(t)$ is bounded, i.e., $\underline \alpha \leq \| x(t) \| \leq \bar \alpha, \forall t\geq 0$
for some constants  $\underline\alpha, \bar\alpha>0$. 
Let $M$ be selected such that $ A_{\zeta} $ in \eqref{M} is Hurwitz
and $\gamma_\zeta$ defined as \eqref{gammazeta}. 
Pick any \begin{align} \label{gammapf}
  0< \gamma <\frac{1}{  \|A\| \|M\|  (\|L U\|  \gamma_\zeta+1 ) }. 
     \end{align}
If there exist matrices $K_o$ and  $P=P\t >0$ satisfying \eqref{PSSP} for 
$\gamma_o  = \min\{\gamma/8 , \gamma/\sqrt{8n}\}$,  
then, for any constants $b_o, b_\zeta >0$,  
 there exist control functions $\mu$ and $\kappa$ such that 
 synchronization is achieve in the sense of \eqref{target} for any initial state of the closed-loop system
 satisfying
   \begin{align}
  \|(W \otimes I_p ) \sigma(0)\| &\leq b_\zeta \nonumber \\
\| \hat \sigma_j^i(0) -\sigma_j(0)\|  &\leq b_o, \; j\in \N_i,\; i \in  \mathcal{V}. \label{ini}
\end{align}

\etheorem

\bproof For the matrix $Q$ in \eqref{lya.equ}, we first define two quantities 
     \begin{align*}
\pi^a_i = &  \|L_i U\|  \| M \|  \sqrt {\varpi(Q)} b_\zeta  \\
& + ( \|L_i U\|  \| M \| \|A \bone\|    \gamma_\zeta  + \| M \|    \|A_i \bone \| ) \\
&\times (\sqrt{ \varpi(P) } b_o +b_x b_K) >0  \end{align*} 
and
     \begin{align*}
\pi^b_i = 1-
 ( \|L_i U\|  \| M \| \|A \bone\|    \gamma_\zeta  + \| M \|    \|A_i \bone \| )  \sqrt{8} \|PB\| \sqrt{ \varpi(P) } 
\end{align*} 
where $b_x$ and $b_K$ are defined in \eqref{bx} and \eqref{bK}, 
$L_i$ and $A_i$ are the $i$-th rows of $L$ and $A$, respectively.

By \eqref{PSSP}, one has 
 \begin{align}
\sqrt{8n} \|PB\| \sqrt{\varpi(P)}&< \gamma .\label{PBpf}
 \end{align} 
 Noting $ \|A \bone\|  \leq\sqrt{n} \|A\|  $ and using \eqref{gammapf} and \eqref{PBpf}, one has
   \begin{align*}
  & ( \|L_i U\|  \| M \| \|A \bone\|    \gamma_\zeta  + \| M \|    \|A_i \bone \| )  \sqrt{8} \|PB\|  \sqrt{ \varpi(P) }\\
<&  \|A \bone\|   \|  M \|  (   \| LU\| \gamma_\zeta  +1   )    \sqrt{8} \|PB\| \sqrt{ \varpi(P) } \\
< & \|A\|    \|  M \|  (   \| LU\| \gamma_\zeta  +1   )    \sqrt{8n} \|PB\|  \sqrt{ \varpi(P) } <1
\end{align*} 
and hence $\pi^b_i  >0$. Now, it is ready to pick    
     \begin{align*}
b_\chi > \frac{\pi^a_i}{\pi^b_i } >0.
\end{align*} 
and hence
 \begin{align*}
 b_\sigma = \sqrt{ \varpi(P) b_o^2+ 8\|PB\|^2 b_\chi^2 \varpi(P) } +b_x b_K.  \end{align*}

Next, we will  prove that the signals $\phi(t)$ and $\chi_i(t)$ are bounded. 
First, it is noted that
 \begin{align*}
 \|\tilde \sigma_j^i(0)\| \leq b_o < b_\sigma. 
 \end{align*}
Assume there exists a finite $T>0$ such that 
 \begin{align*}
 \|\tilde \sigma_j^i(t)\|  <  b_\sigma,\; \forall t <T
 \end{align*}
 and
  \begin{align}
  \|\tilde \sigma_j^i(T)\|  =  b_\sigma. \label{hatsigmaT}
 \end{align}
 A contradiction is derived below.

By the assumption, one has \begin{align}
   \|\phi_i(t)  \| & \leq    \sum_{j \in \N_i} a_{ij}  \| \tilde\sigma^i_j(t)\|    \leq   \sum_{j \in \N_i} a_{ij}  b_\sigma 
   =  \|A_i \bone \| b_\sigma
  \end{align}
and hence 
\begin{align}
  \|\phi(t)  \| & \leq  \|A \bone\| b_\sigma,
  \end{align}
  which shows that $\phi(t)$ is bounded for   $t\in [0,T]$.
   
 Using \eqref{Vzeta} and \eqref{dVzeta} gives
      \begin{align*} \dot V(\zeta ) \leq  - \frac{V(\zeta )}{ \lambda_{\max}(Q) } + \| Q (W  \otimes BM)\|^2  \| \phi\|^2. 
\end{align*}
As a result, for   $t\in [0,T]$,
   \begin{align*}
V(\zeta(t)) & \leq   V(\zeta (0) ) +\| Q (W  \otimes BM)\|^2 \|A \bone\|^2 b_\sigma^2  \lambda_{\max}(Q) 
\end{align*}
and
   \begin{align*}
 \|\zeta(t)\|^2  \leq  &
 \frac{\lambda_{\max}(Q)  \|\zeta (0)\|^2  } {\lambda_{\min}(Q) }  \\ & +\| Q (W  \otimes BM)\|^2 \|A \bone\|^2 b_\sigma^2 \frac{\lambda_{\max}(Q) }{\lambda_{\min}(Q) }\\
  \leq  & \varpi(Q) \|\zeta (0)\|^2      +  \|A \bone\|^2 b_\sigma^2   \gamma_\zeta^2.
\end{align*}

Then, we use the facts
  \begin{align*}
  (L_i\otimes I_p ) \sigma = (L_i T^{-1} \otimes I_p ) \left[\begin{array}{c}
\bar \sigma\\
\zeta
\end{array}\right] =(L_i U \otimes I_p )\zeta
\end{align*}
and
     \begin{align}
\chi_i  
= -(L_i \otimes M)  \sigma    + M \phi_i,
\end{align}
and have the following calculation, for   $t\in [0,T]$,
   \begin{align*}
\|\chi_i (t)\|  
\leq & \| (L_i \otimes M)  \sigma(t) \|   + \| M \|  \| \phi_i (t)\| \\
\leq & \| L_i U \otimes M \| \|\zeta(t) \|   + \| M \|  \| \phi_i (t)\|  \\
\leq & \|L_i U\|  \| M \|  \sqrt {\varpi(Q) \|\zeta (0)\|^2     +  \|A \bone\|^2 b_\sigma^2   \gamma_\zeta^2} \\  
& + \| M \|    \|A_i \bone \| b_\sigma  \\
  \leq  & \|L_i U\|  \| M \|  \sqrt {\varpi(Q)}  b_\zeta   \\
&+ ( \|L_i U\|  \| M \|
\|A \bone\|    \gamma_\zeta  + \| M \|    \|A_i \bone \| ) b_\sigma 
\end{align*} 
and hence
   \begin{align*}
\|\chi_i (t)\|  
\leq  & \|L_i U\|  \| M \|  \sqrt {\varpi(Q)} b_\zeta   \\
& + ( \|L_i U\|  \| M \|
\|A \bone\|    \gamma_\zeta  + \| M \|    \|A_i \bone \| )  \\
& (\sqrt{ \varpi(P) b_o^2+8\|PB\|^2  b_\chi^2 \varpi(P)  } +b_x b_K) \\
\leq & \|L_i U\|  \| M \|  \sqrt {\varpi(Q)} b_\zeta   \\
& + ( \|L_i U\|  \| M \|
\|A \bone\|    \gamma_\zeta  + \| M \|    \|A_i \bone \| )  \\
& (\sqrt{ \varpi(P) } b_o +  \sqrt{8} \|PB\| \sqrt{ \varpi(P) }  b_\chi   +b_x b_K)\\
=&  \pi^a_i  + (1-\pi_i^b) b_\chi< b_\chi.
\end{align*} 

Denote $\bar b_\chi = \max_{i=1,\cdots, n} \{ \pi^a_i  + (1-\pi_i^b) b_\chi \} < b_\chi$, 
one has
   \begin{align*}
 \|\chi_i (t)\|  \leq \bar b_\chi, \; \forall t\in [0,T].
\end{align*} 
 Using the similar development of \eqref{nubound2}  gives
  \begin{align}
 \|\nu(t)\|^2 \leq \varpi(P) b_o^2+ 8\|PB\|^2 \bar b_\chi^2 \varpi(P). \end{align}
Moreover,  using \eqref{bK} and \eqref{xbound} gives
 \begin{align*}
 \|\tilde \sigma_j^i(t)\| \leq &
 \|\nu(t)\| +  \|K(\hat x(t))  {\tilde x(t)}\| \\
 \leq & \sqrt{ \varpi(P) b_o^2+ 8\|PB\|^2 \bar b_\chi^2 \varpi(P) } +b_x b_K
  < b_\sigma
 \end{align*}
for $ t\in [0,T]$. It contradicts with the assumption \eqref{hatsigmaT}. 
Using the proof by contradiction,  one has 
 \begin{align}  \|\tilde \sigma_j^i(t)\|  <  b_\sigma,\;
  \|\phi(t)  \|  \leq  \|A \bone\| b_\sigma, \;
  \|\chi_i (t)\|  \leq b_\chi ,\; \forall t>0.
  \end{align}

Finally, Lemmas~\ref{lemma:sigma} and \ref{lemma:chi} can be applied to complete the proof
using a small gain argument. 
On one hand, by Lemma~\ref{lemma:sigma},
\begin{align*}
\limsup_{t\rightarrow \infty} \| \tilde \sigma_j^i(t)  \| & \leq \gamma \limsup_{t\rightarrow \infty} \|\chi_j(t)\| ,\;
 j \in \N_i.
\end{align*}
For  $\phi_i = \sum_{j \in \N_i} a_{ij}   \tilde\sigma^i_j$,  one has
\begin{align*} 
\limsup_{t\rightarrow \infty}  \|\phi_i(t)  \| & \leq \gamma  \sum_{j \in \N_i} a_{ij}   \limsup_{t\rightarrow \infty} \|\chi_j(t)\|  
  \end{align*}
and hence
\begin{align}\label{sgpf1}
 \limsup_{t\rightarrow \infty}  \|\phi(t)  \| & \leq \gamma \|A\|  \limsup_{t\rightarrow \infty} \|\chi(t)\|. 
  \end{align}
 On the other hand, by Lemma~\ref{lemma:chi}, 
    \begin{align}\label{sgpf2}
\limsup_{t\rightarrow \infty}  \|\chi (t) \|
\leq   \gamma_\chi \limsup_{t\rightarrow \infty}  \| \phi(t)\|
\end{align} 
for 
$\gamma_\chi  = \|M\| (\| L U \| \gamma_\zeta +1)$ that satisfies the small gain condition $ \gamma \|A\| \gamma_\chi <1$. 
Putting \eqref{sgpf1} and \eqref{sgpf2} together gives 
\begin{align} 
 \limsup_{t\rightarrow \infty}  \|\phi(t)  \| & \leq \gamma \|A\|    \gamma_\chi \limsup_{t\rightarrow \infty}  \| \phi(t)\|.
  \end{align}
It concludes that \begin{align*} 
 \lim_{t\rightarrow \infty}  \|\phi(t)  \| =0.
  \end{align*}
 Using the same approach, one has 
\begin{align*}  \lim_{t\rightarrow \infty}  \|\chi (t) \|=0.\end{align*} 
 Finally, by \eqref{limzeta}, one has
 \begin{align*}
\lim_{t\rightarrow \infty}\| \zeta(t)\|=0
\end{align*}
 and
 \begin{align*}
\lim_{t\rightarrow \infty} \| \sigma(t) - (\bone \otimes I_p)  \bar\sigma(t)\| = \lim_{t\rightarrow \infty} 
\| (U \otimes I_p) \zeta(t) \|=0,
\end{align*}
which implies \eqref{target}. 
 The proof is thus completed. \eproof
  
\bremark  The functions $\mu$ and $\kappa$ used in the observer
\eqref{observeri} can be explicitly designed accordingly Lemma~\ref{lemma:sigma}. 
From the proof, it is noted that the functions depend on the function $f$ and the parameters
$S$, $E$, $B$ $\alpha$, and $b_\chi$. Also, it is noted from the proof of 
Theorem~\ref{thm-FMMAS} that  $b_\chi$ further depends on $L$, $b_\zeta$, and $b_o$. 
Therefore, the two functions $\mu$ and $\kappa$ can be uniform for all the agents as shown in \eqref{observeri}. 
As the constants $b_\zeta$ and $b_o$ to characterize of the initial states of the closed-loop system can be arbitrarily 
selected, the synchronization problem in the theorem is solved in a semi-global sense. 
\eremark

\bremark The state of the dynamics \eqref{MAS} of agent $i$ is $(\sigma_i, x_i)$, 
where $x_i$ is accessible by the neighbored agents via network but $\sigma_i$ is not. 
Synchronization of the agents in terms of the full state trajectories $(\sigma_i(t), x_i(t))$ obviously implies 
frequency synchronization \eqref{target} for $\sigma_i$. It is worth mentioning that 
the mechanism of frequency synchronization developed in this paper does not follow the traditional idea of 
synchronization of full state trajectories used in most existing works.  
In particular, the difference between the trajectories of two agents does not approach zero, that is,
$\lim_{t\rightarrow\infty} \| x_i(t) -x_j(t) \|$ does not necessarily exist or is not equal to zero, even when  \eqref{target} is guaranteed. 
  \eremark

\section{Numerical Examples} \label{section:example}

\begin{figure}[h]
\centering
    \includegraphics[width=0.5\hsize]{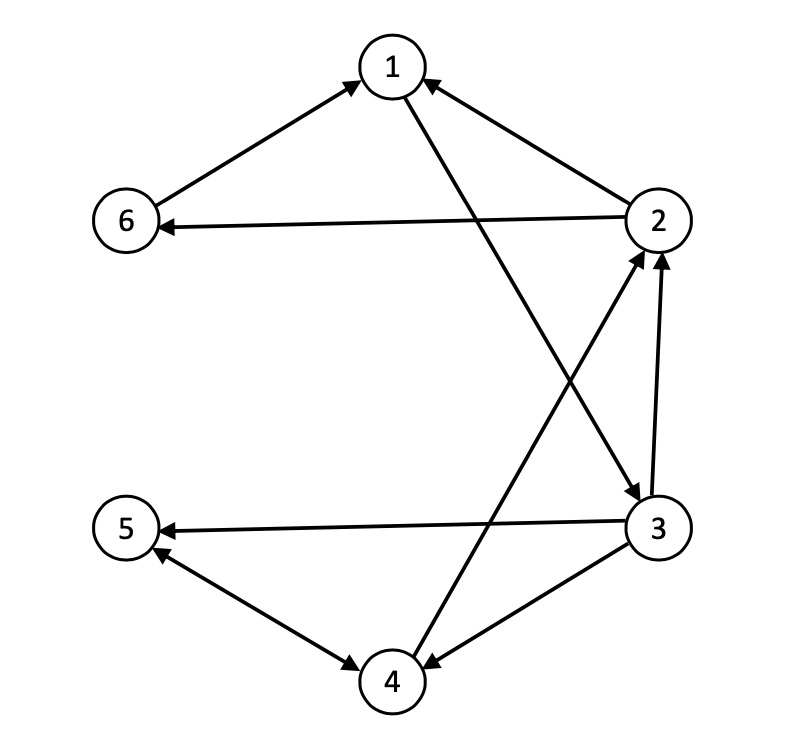}
  \caption{Illustration of the network topology of six oscillators (agents) where the weight associated with
  each edge (arrowed line) is positive. }\label{fig.topo}
\end{figure}

We consider a network of $n=6$ frequency modulated oscillators of the form \eqref{MAS}
with 
 \begin{align*}
  S=\left[\begin{array}{cc} 0 & 0.1 \\-0.1 & 0 \end{array} \right], \;
  B=\left[\begin{array}{c}  1 \\ 1 \end{array} \right],\;
    E=\left[\begin{array}{cc}  4.5 &  0 \end{array} \right] 
 \end{align*}
 and
 \begin{align*}
 f(x_i)=  
  \left[\begin{array}{cc}   0 & 1  \\-1 &0   \end{array} \right] x_i,\; f_o(x_i) =0, \;
  \omega_c=3.
   \end{align*} 
 In this example, the modulated signals $x_i(t)$ are frequency-varying sinusoids.
 The target is to synchronize the frequencies of the signals generated by the oscillators, not the instantaneous  values of the signals.  
 The network topology is illustrated in Fig.~\ref{fig.topo}
   
In particular, the following function 
     \begin{align*}
 K(\hat x)= K_o \frac{ \hat x\t }{\| \hat x \|^2 }
  \left[\begin{array}{cc}   0 &-1  \\1 &0   \end{array} \right],\;  K_o = \left[\begin{array}{c}  1.80\\ 1.76  \end{array} \right]    \end{align*}
with   
 \begin{align*}
K'(x, \hat x, \hat\sigma)= & K_o \left(  \frac{   I_q
- 2 f(\hat x)  f\t(\hat x)   
 } {\| \hat x \|^2 }  \times \right. \\
 & \left. \left[\begin{array}{cc}   0 & 1  \\-1 &0   \end{array} \right]   
 \left(f(x) \hat \omega    +\kappa (x , \hat x ) \right)
 \right)\t 
   \end{align*}
is used to construct the functions $\mu$ and $\kappa$ of the frequency observers. 
The parameter $\beta=10$ is used for the simulation. It is noted that a higher value of $\beta$ gives faster 
convergence of the frequency observers.  For $a_{ij}=1,\; \forall \left(i, j\right) \in \mathcal{E}$, 
the control gain matrix $M=0.01  \left[\begin{array}{cc} 1 & 0.5\end{array} \right]$ representing 
weak coupling among the oscillators is used in the controller. 
Using the observer/controller proposed in the paper, the synchronization performance of the oscillators 
is discussed below. In the figures, we only demonstrate the trajectories of three oscillators, denoted
as OSC 1, OSC 3, and OSC 5, corresponding to the nodes 1, 3, and 5 in Fig.~\ref{fig.topo}, respectively.
The other three oscillators have similar profiles and thus omitted.

\begin{figure}[t]
\centering
    \includegraphics[width=1\hsize]{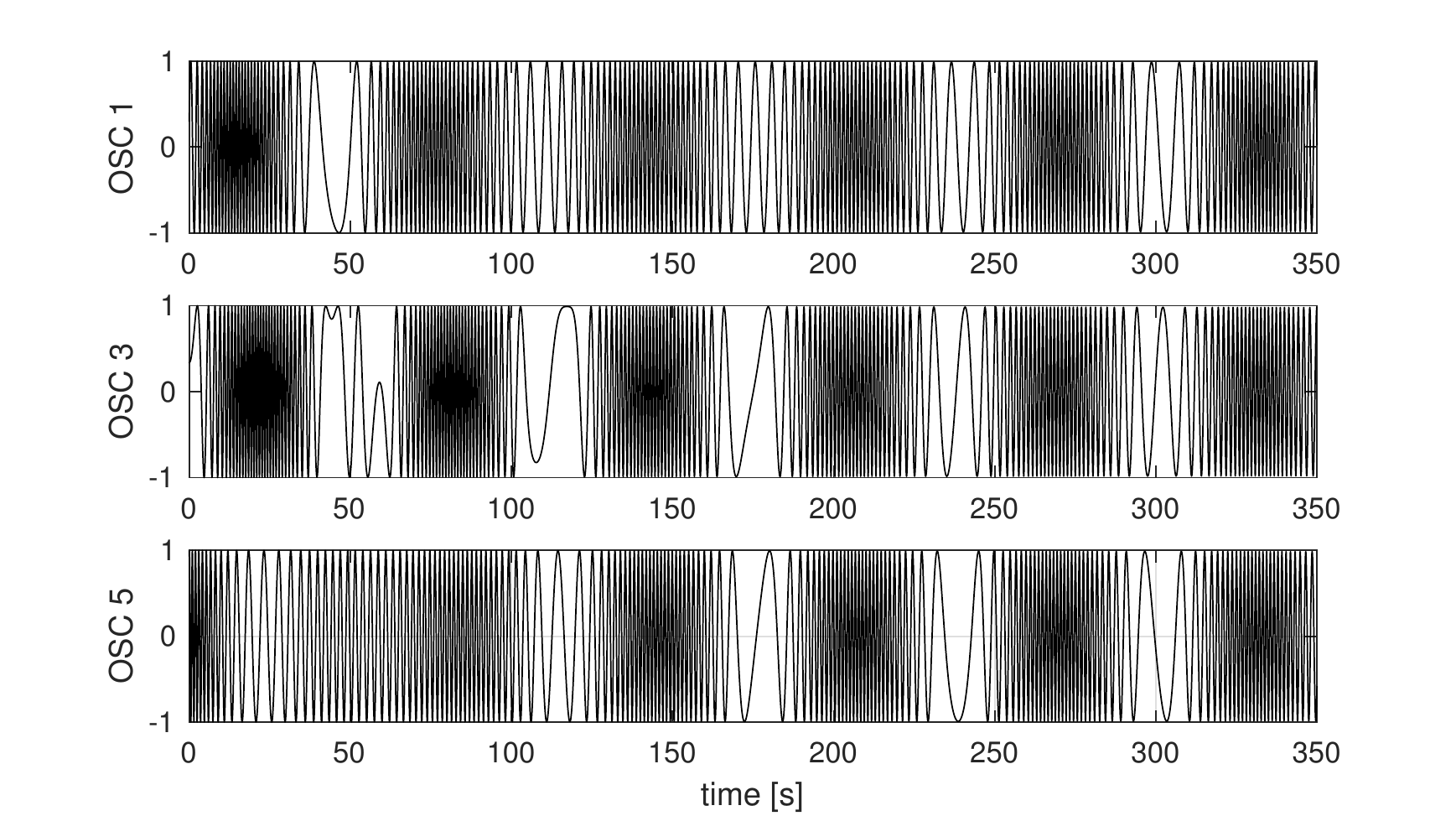}
  \caption{Profile of the frequency modulated signals $x_{i1}(t)$, the first element of the state vector $x_i(t)$,
  for $i=1,3,5$.}\label{fig.exmp1x}
%
\centering
    \includegraphics[width=1\hsize]{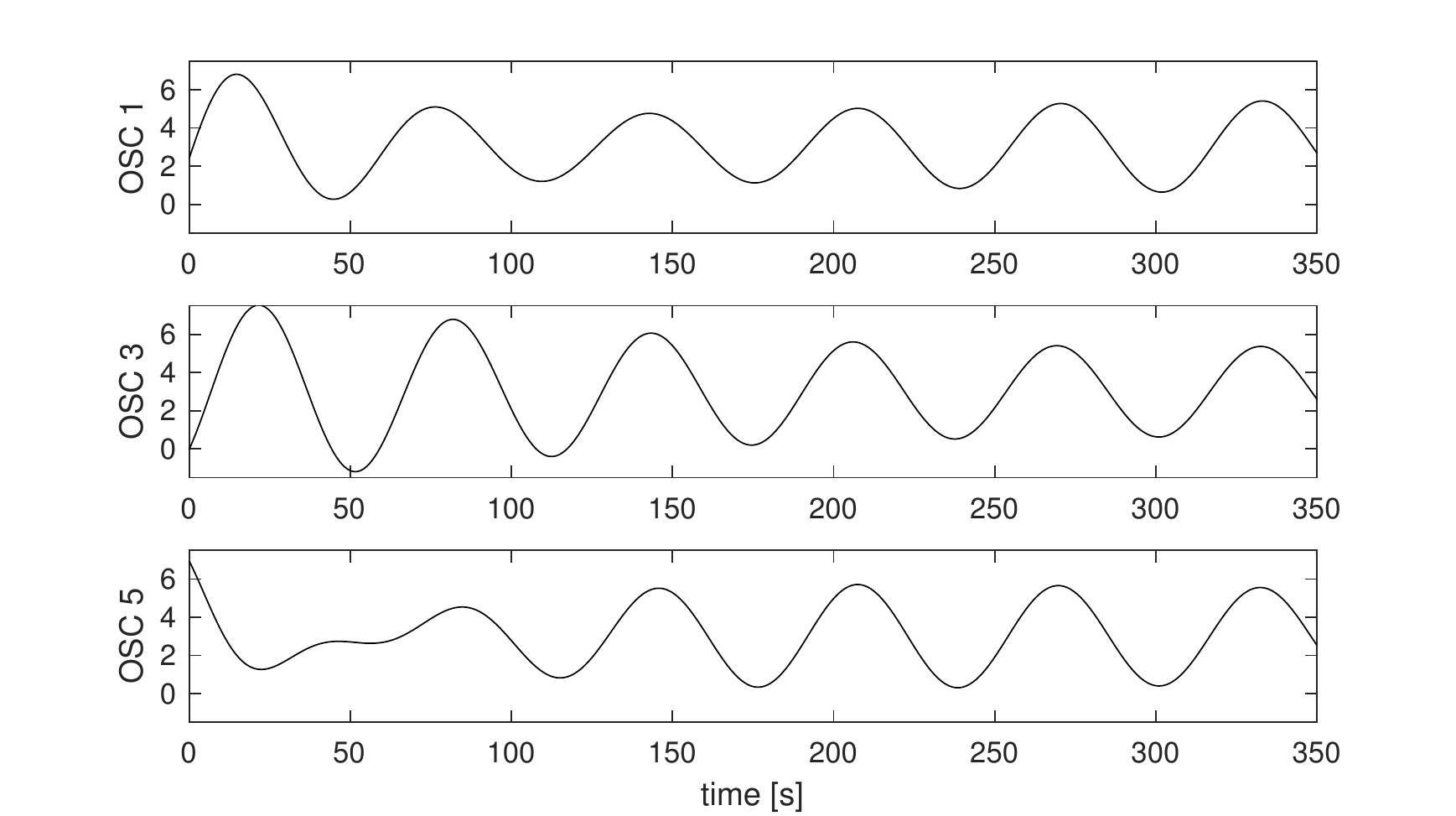}
  \caption{Profile of the frequencies $\omega_i(t)$ for $i=1,3,5$.}
\label{fig.exmp1sigma}
%
\centering
    \includegraphics[width=1\hsize]{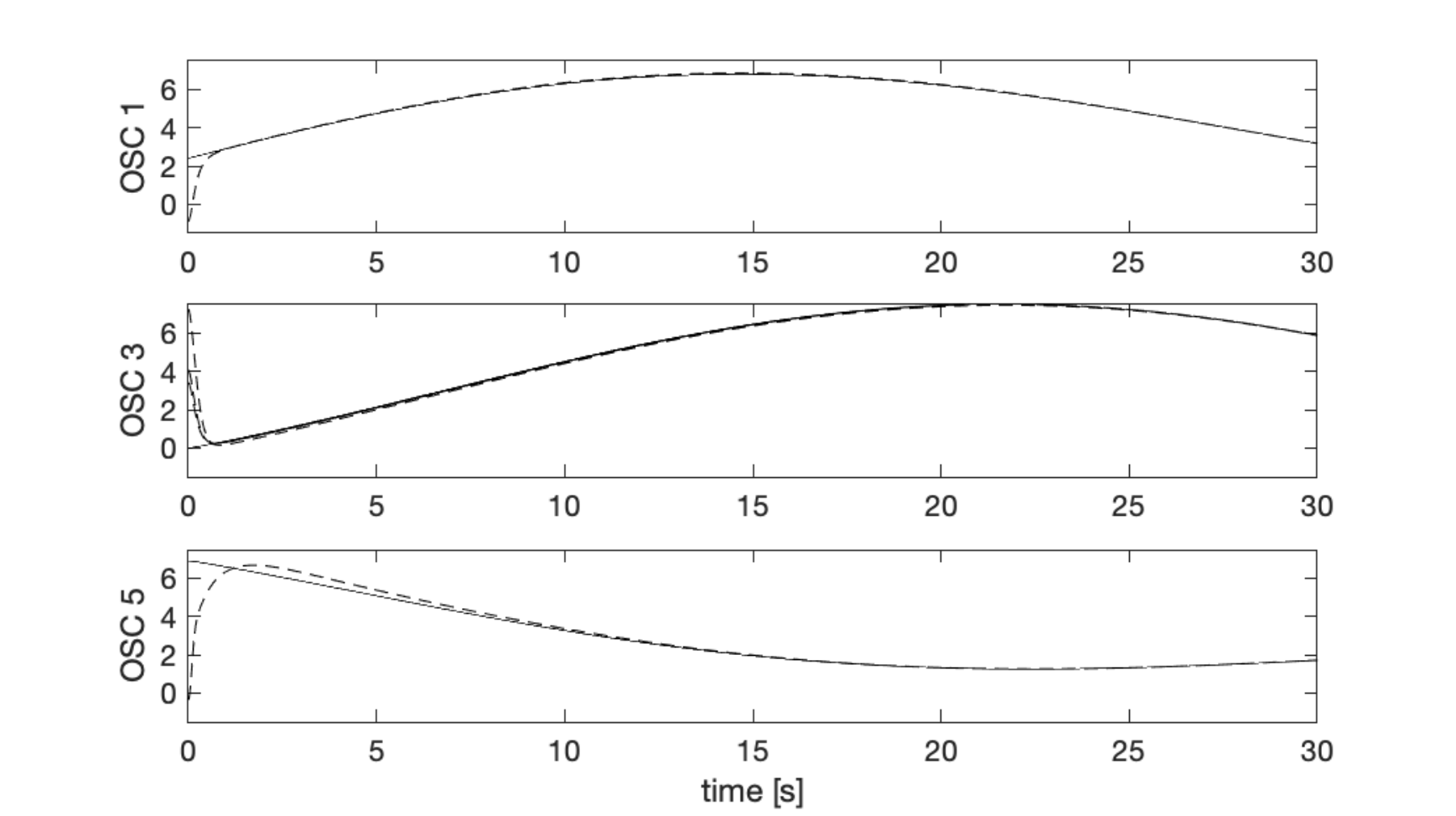}
  \caption{Profile of the observed frequencies $\hat\omega_j^i(t)$, $j \in \N_i$ (dashed curves) vs the real frequencies $\omega_j(t)$ (solid curves), for $j=1,3,5$.  }
  \label{fig.exmp1error}
\end{figure}

The profile of the frequency modulated signals $x_i(t)$ (only the first 
component $x_{i1}(t)$ for concise presentation) is shown in Fig.~\ref{fig.exmp1x}. 
Frequency synchronization can be observed in this figure as expected, while synchronization
of the signals in oscillation phase and amplitudes is not needed. 
Frequency synchronization can also be exhibited in Fig.~\ref{fig.exmp1sigma}
in an explicit illustration. It is noted that frequencies $\omega_i(t)$ in this figure are 
not directly measured or transmitted via the network. 
The effectiveness of the proposed frequency observer is verified in 
Fig.~\ref{fig.exmp1error}. 
According to the topology, the frequency of OSC~1 is observed by agent 3,
that of OSC~3 is observed by agents 2, 4, 5, separately,  and that of OSC~5 by agent 4. 
The convergence of the observed frequencies (in dashed curves) to the real values (solid curves) 
can be achieved by all the observers . 

\begin{figure}[t]
\centering
    \includegraphics[width=1\hsize]{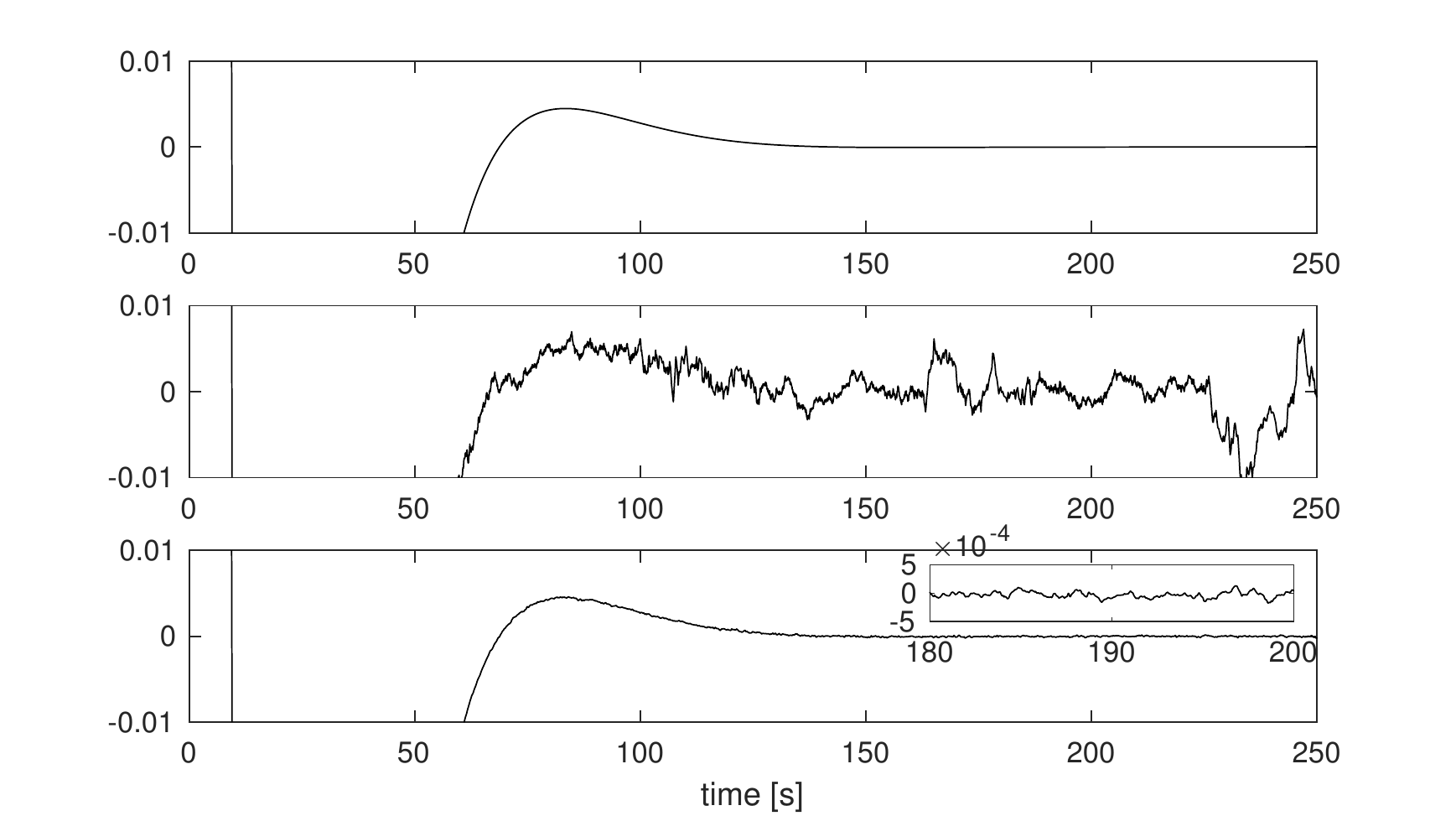}
  \caption{Profile of the synchronization error $\omega_1(t) -\omega_3(t)$ in three cases.
 Top graph: noise free; middle graph: unmodulated synchronization with noise; and bottom graph:
 frequency modulated synchronization with noise.}\label{fig.exmp1noise}
\end{figure}

Next, we examine the robustness of synchronization of frequency modulated oscillators with respect to 
noise in the transmitted signals via the network.
The result is shown in Fig.~\ref{fig.exmp1noise} using  
the synchronization error  $\omega_1(t) -\omega_3(t)$ as an example.  
The top graph shows a perfect synchronization behavior without any noticeable error in the scale of $\pm 0.01$.
The relative and absolute error tolerances in MATLAB simulation are $10^{-4}$ and $10^{-8}$, respectively. 
The middle graph shows the result for the unmodulated controller \eqref{chiideal} with transmission noise, 
 i.e, 
 \begin{align*}  
\chi_i = -  M \sum_{j \in \N_i}  a_{ij} (\sigma_i -\sigma_j-\n^i_j)  ,\; i \in  \mathcal{V} 
\end{align*}
where the level of noise $\n^i_j$ on $\sigma_j$ transmitted from $j$ to $i$ is $\pm 1\%$ of the magnitude of the transmitted signals.
In this case, the synchronization error induced by the noise can ben obviously observed. 
 Finally, we apply the same level of noise on the frequency modulated controller \eqref{observeri}, i.e., 
\begin{align*}
\dot {\hat \sigma}_j^i &= S \hat \sigma_j^i +\mu (x_j+\n^i_j, \hat x_j^i, \hat\sigma_j^i) \nonumber\\
\hat\omega_j^i &= E\hat\sigma_j^i +\omega_c  \nonumber\\
\dot {\hat x}_j^i&=   f( \hat x^i_j)\hat \omega_j^i +f_o(x_j+\n^i_j)  +\kappa (x_j+\n^i_j, \hat x_j^i ),\;
j\in \N_i .
\end{align*}
The result in the bottom graph exhibits  that the synchronization error using frequency modulation 
is significantly less than that in  the unmodulated case. 
The error is noticeable with more details when it is amplified to the scale of $10^{-4}$. 
The comparison in Fig.~\ref{fig.exmp1noise} verifies the robustness of frequency modulated MASs
in mitigating the influence of network noise.

 \medskip
   
Next, we consider another example of a network of Hindmarsh-Rose oscillators, which are a class of
low-dimensional  Hodgkin-Huxley model of neuronal activity. They can characterize  the spiking-bursting behavior of the membrane potential observed in experiments. The dynamics are of the form \eqref{MAS} with
 \begin{align*}
  S=\left[\begin{array}{cc} 0 & 0.005 \pi \\-0.005 \pi & 0 \end{array} \right], \;
  B=\left[\begin{array}{c}  1 \\ 1 \end{array} \right],\;
    E=\left[\begin{array}{cc}  0.4 &  0 \end{array} \right] 
 \end{align*}
 and
 \begin{align*}
x_i=  
  \left[\begin{array}{c}
v_i \\
w_i\\
z_i\end{array} \right],\;  f(x_i)=  
  \left[\begin{array}{c}
2\\
-5 v_i^2 - w_i +1\\
0
\end{array} \right],\\
f_o(x_i)=  
  \left[\begin{array}{c}
3 v_i^2-v_i^3+w_i-z_i \\
0\\
 0.005(4(v_i+1.5)-z_i)
\end{array} \right],\; \omega_c=0.9.
   \end{align*} 
The variable $v_i(t)$ represents the membrane potential of a neuron,
 $w_i(t)$ is  the spiking variable that represents the transport rate of sodium and potassium ions through fast ion channels,
 and $z_i(t)$ is the adaptation current that increases at every spike and leads to a decrease in the firing rate.
The signal $\omega_i(t)$ is encoded as the firing rate of the membrane potential. 
 The network topology is the same one in Fig.~\ref{fig.topo}.
 The non-sinusoidal signals $x_i(t)$ of  this class of oscillators
 consist of bursts and spikes of varying  firing rates. 
 The synchronization target is not the instantaneous signals but their firing rates. 
 
\begin{figure}[t]
\centering
    \includegraphics[width=1\hsize]{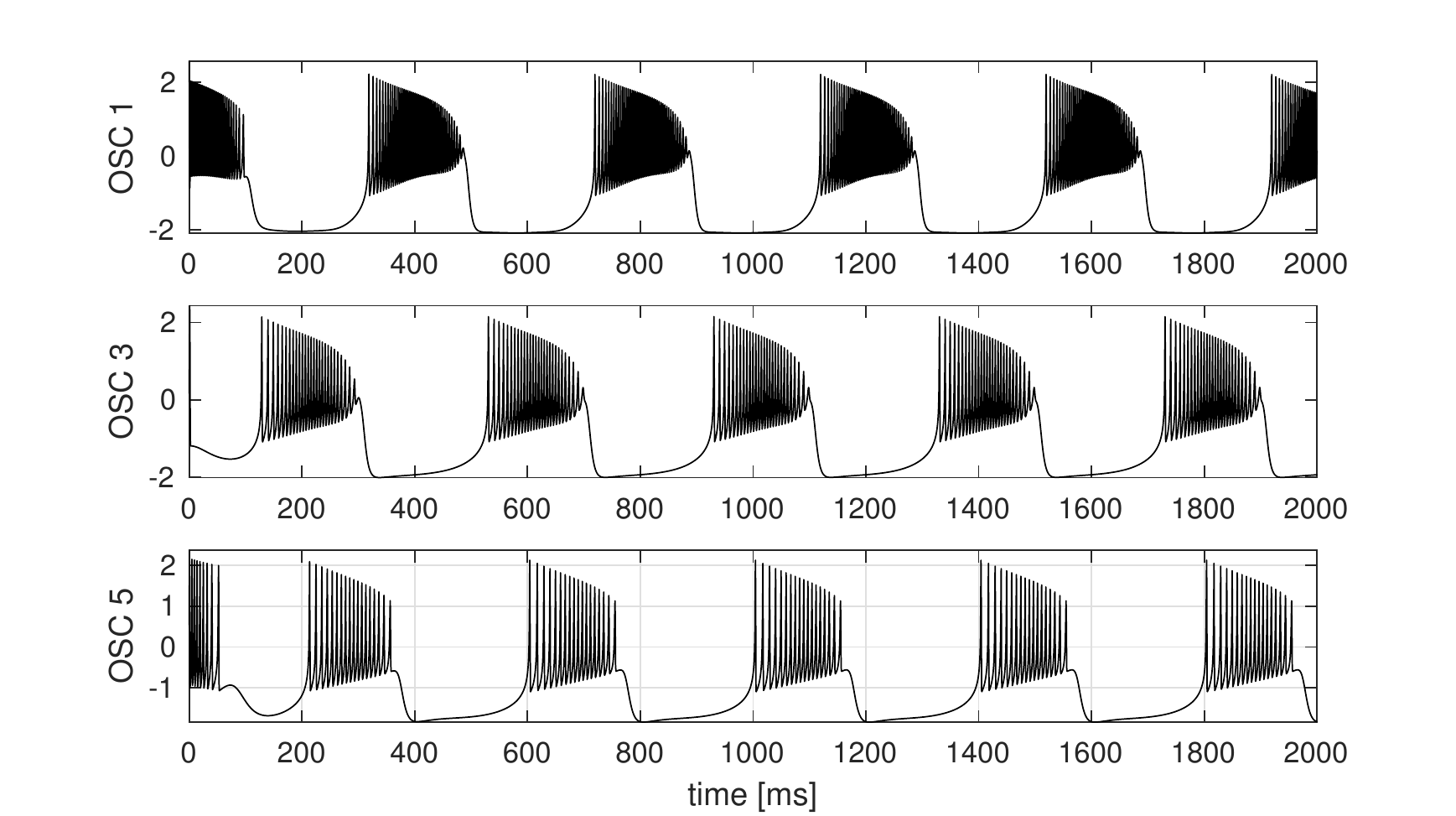}
  \caption{Profile of the membrane potentials $v_{i}(t)$ of individual Hindmarsh-Rose oscillators  without synchronization for $i=1,3,5$.}\label{fig.exmp2nocoupling}
\end{figure}

\begin{figure}[t]
\centering
    \includegraphics[width=1\hsize]{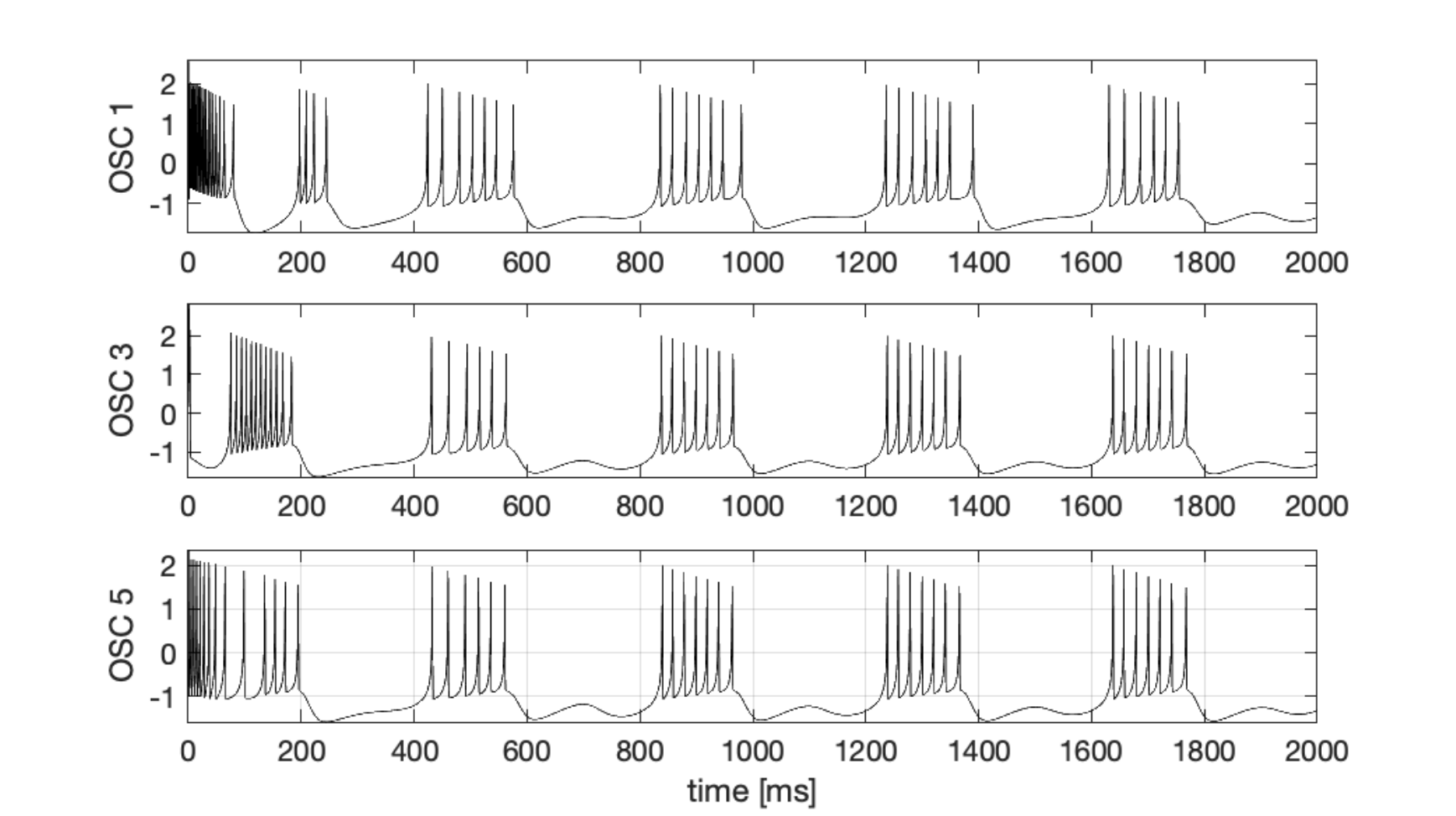}
  \caption{Profile of the membrane potentials $v_{i}(t)$ of Hindmarsh-Rose oscillators  for $i=1,3,5$.}\label{fig.exmp2x}
%
\centering
    \includegraphics[width=1\hsize]{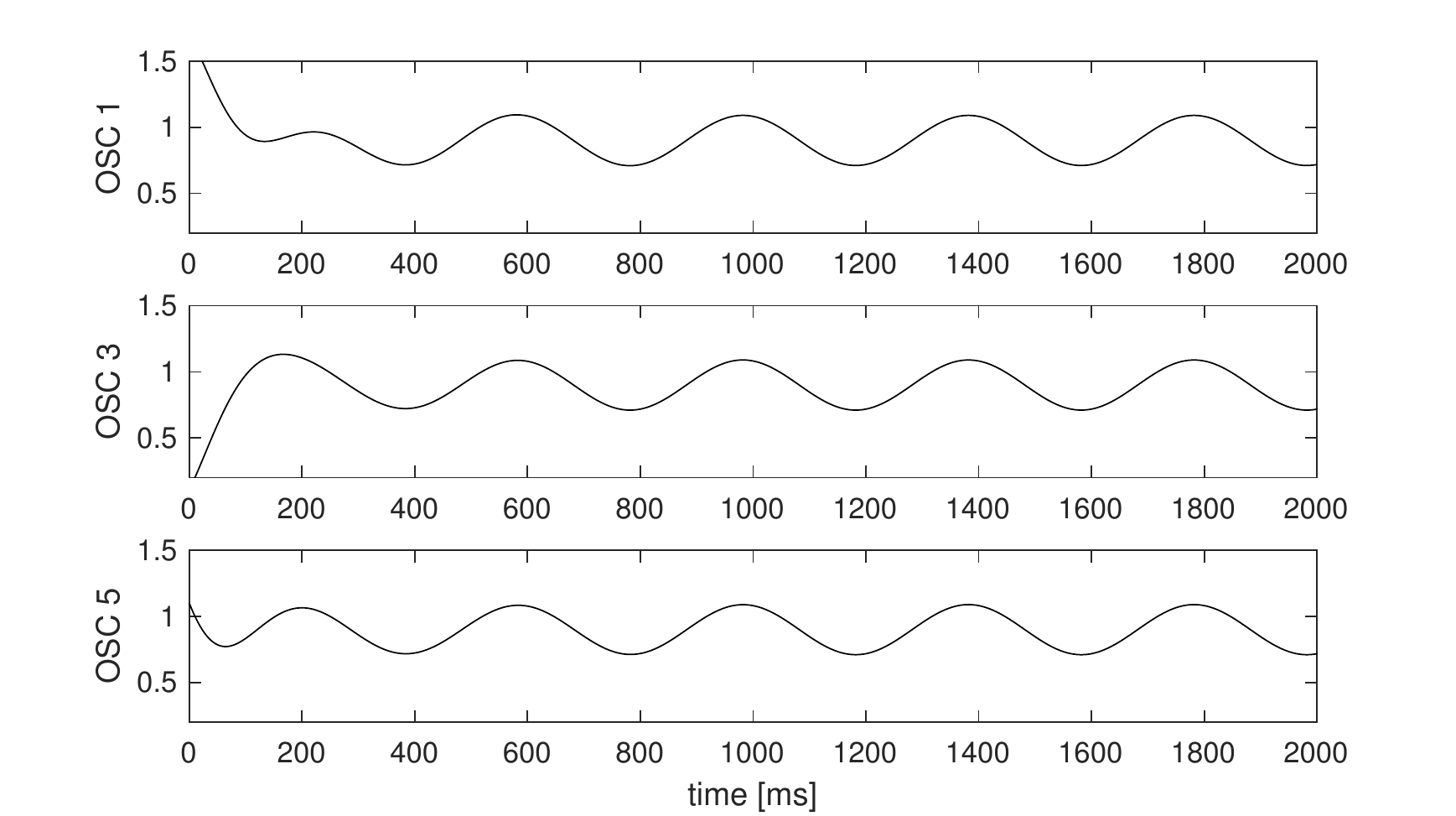}
  \caption{Profile of the firing rates $\omega_i(t)$ of Hindmarsh-Rose oscillators  for $i=1,3,5$.}\label{fig.exmp2sigma}
    \centering
    \includegraphics[width=1\hsize]{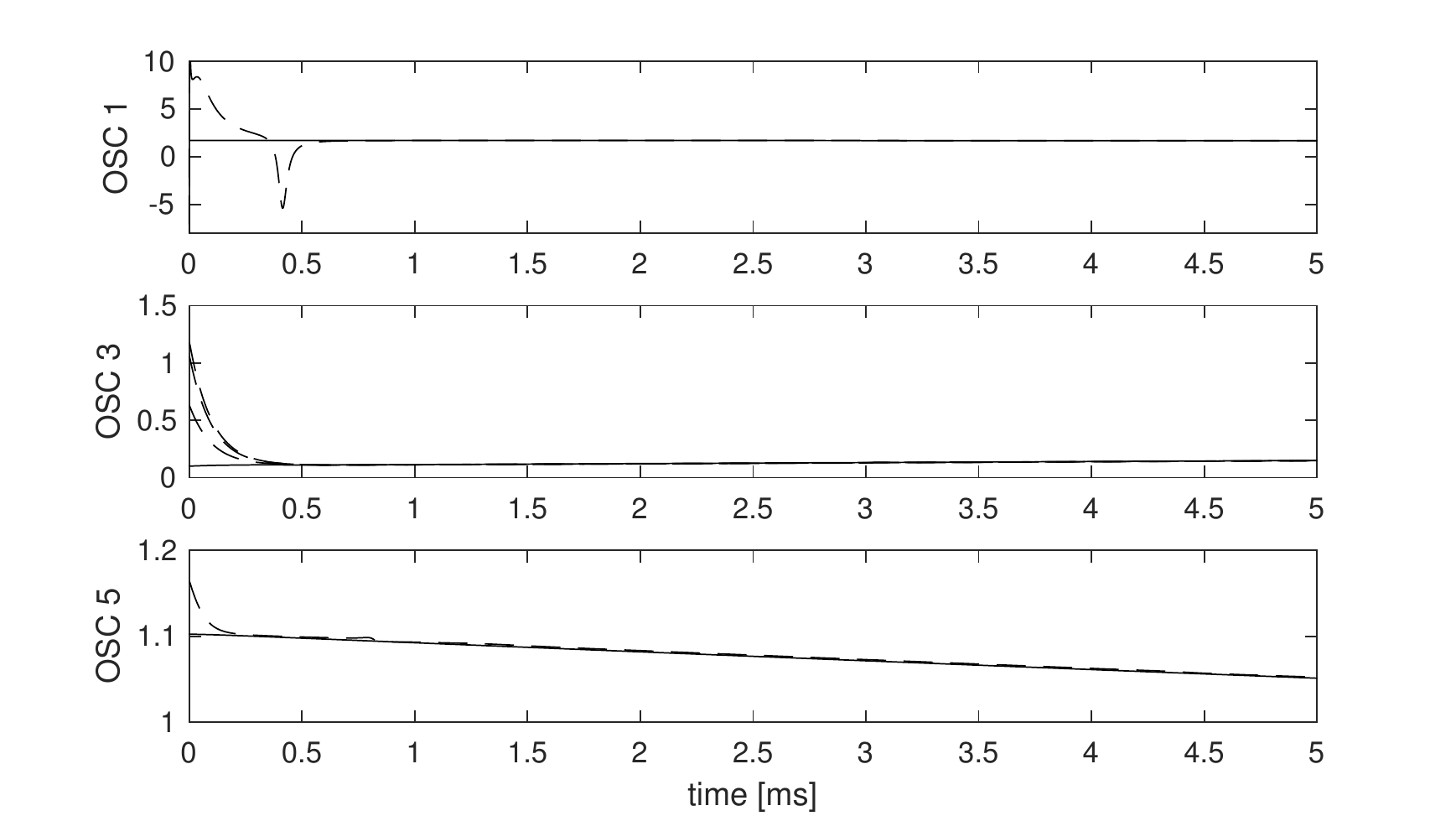}
  \caption{Profile of the observed  firing rates $\hat\omega_j^i(t)$, $j \in \N_i$ (dashed curves) vs the real  firing rates $\omega_j(t)$ (solid curves) of Hindmarsh-Rose oscillators for $j=1,3,5$.  }
  \label{fig.exmp2error}
\end{figure}

The results discussed in the first example can be observed in this example as
well.  As a comparison,  the membrane potentials of 
 the individual oscillators without synchronization control, i.e., $\chi_i=0$,  are plotted in Fig.~\ref{fig.exmp2nocoupling}.
With the synchronization controller proposed in this paper,
frequency synchronization of membrane potentials  is plotted in Fig.~\ref{fig.exmp2x}
and synchronization of firing rates of spikes in Fig.~\ref{fig.exmp2sigma}.  
Once synchronized, the firing rates of spikes are generally less than those of individual oscillators 
because the variations of the firing rates are averaged. 
It is also noted that the signals in terms of  the instantaneous values are not necessarily synchronized. 
The convergence of the observed firing rates to the real values is shown in  Fig.~\ref{fig.exmp2error}.
The results well demonstrate
the similar phenomenon observed in real neural circuits. 
 
%
%
%
%
%
%
%
%
%
%
%
%
%
%
%
%
%
%
%

\section{Conclusions}\label{section:conclusion}

In this paper, we have established a framework for synchronization of frequency modulated multi-agent systems for the first time.
The framework includes a frequency observer subject to an external input and a distributed  synchronization controller based
on the observer using the small gain argument. We have  verified the effectiveness of the design using rigorous theoretical proofs and extensive numerical simulation. 
It is interesting to further extend the framework to study more complicated models, especially from real biological neural networks. 
The research for unmodulated multi-agent systems has experienced rapid development in the recent decade, which can be further developed for frequency modulated  multi-agent systems using the framework in this paper. 
The robust synchronization mechanism using frequency modulated signals has the potentials in 
revealing the similar phenomenon in real systems and gaining deep understanding of how 
neuronal circuits generate extremely robust and adaptive oscillatory behaviors in animal locomotion. 
It is expected that the research along this direction is also beneficial for
building artificial frequency modulated multi-agent systems with more engineering applications.

\bibliographystyle{IEEEtran}

\end{document}